\documentclass[a4paper,12pt]{article}
\usepackage{amsmath,amsthm,amssymb,latexsym,epsfig,graphicx}
\title{{\bf Extra cancellation of even Calder\'{o}n-Zygmund operators and Quasiconformal mappings }}
\author{\Large{\Large Joan Mateu, Joan Orobitg and Joan Verdera}}
\newtheorem{teorema}{Theorem}
\newtheorem*{teor}{Theorem}
\newtheorem*{teor'}{Theorem'}

\newtheorem*{cor1}{Corollary 1'}
\newtheorem*{cor2}{Corollary 2'}
\newtheorem{co}{Corollary}
\newtheorem{lemma}[co]{Lemma}
\newtheorem*{ML}{Main Lemma}

\newtheorem*{sublemma}{Sublemma}

\theoremstyle{definition}
\newtheorem*{example}{Example}

\newtheorem*{gracies}{Acknowledgements}

\newcommand{\Rn}{{\mathbb R}^n}
\newcommand{\ep}{\varepsilon}
\newcommand{\C}{\mathbb{C}}

\newcommand{\Lip}{\operatorname{Lip}(\ep',\Omega)}
\newcommand{\Lipo}{\operatorname{Lip}(\ep,\Omega)}
\newcommand{\fr}{\partial\,\Omega}
\newcommand{\R}{r_0}

\begin{document}

\date{}

\maketitle

\begin{abstract}
In this paper we discuss a special class of Beltrami coefficients
whose associated quasiconformal mapping is bilipschitz. A particular
example are those of the form $f(z)\chi_\Omega(z)$, where $\Omega$
is a bounded domain with boundary of class $C^{1+\ep}$ and $f$ a
function in $\operatorname{Lip}(\ep,\Omega)$ satisfying $\|
f\|_\infty < 1$. An important point is that there is no restriction
whatsoever on the $\operatorname{Lip}(\ep,\Omega)$ norm of $f$
besides the requirement on Beltrami coefficients that the supremum
norm be less than $1$. The crucial fact in the proof is the extra
cancellation enjoyed by even homogeneous Calder\'{o}n-Zygmund kernels,
namely that they have zero integral on half the unit ball. This
property is expressed in a particularly suggestive way and is shown
to have far reaching consequences.

An application to a Lipschitz regularity result for solutions of
second order elliptic equations in divergence form in the plane is
presented.
\end{abstract}

\section{Introduction}
Consider the Beltrami equation
\begin{equation}\label{eq1}
\frac{\partial \Phi}{\partial \overline{z}}(z)= \mu(z)\,
\frac{\partial \Phi}{\partial z}(z)\,,\quad z \in \C\,,
\end{equation}
where $\mu$ is a Lebesgue measurable function on the complex plane
$\C$ satisfying $\|\mu\|_\infty < 1$\,. According to a remarkable
old theorem of Morrey \cite{M} there exists an essentially unique
function $\Phi$ in the Sobolev space $W^{1,2}_{\text{loc}}(\C)$
\,(functions with first order derivatives locally in $L^2$) which
satisfies \eqref{eq1} almost everywhere and is a
\emph{homeomorphism of the plane}. These functions are called
quasiconformal. It turns out that $\Phi$ may change drastically
the Hausdorff dimension of sets. Indeed, sets of arbitrarily small
positive Hausdorff dimension may be mapped into sets of Hausdorff
dimension as close to $2$ as desired (and the other way around by
the inverse mapping). There has been during the last decades much
hard and penetrating work in understanding how $\Phi$ distorts
sets (see, for instance \cite{As} and the references given there
or \cite{LSU} for a recent result).

On the other hand, orientation preserving bilipschitz homeomorphisms
of the plane are easily seen to satisfy a Beltrami equation for a
certain Beltrami coefficent $\mu.$   Since bilipschitz mappings
preserve all metric properties of sets, in particular Hausdorff
dimension, they appear to be a distinguished subclass of
particularly simple quasiconformal mappings . In \cite{R} one gives
geometric conditions which are necessary and sufficient for $\Phi$
being bilipschitz, but which do not involve the Beltrami coefficient
$\mu$. In fact, it is widely accepted that the problem of
characterizing in an efficient way those $\mu$ which determine
bilipschitz mappings is hopeless.

A classical result that goes back to Schauder (\cite{AIM})  asserts
that $\Phi$ is of class~$C^{1+\varepsilon}$ provided $\mu$ is a
compactly supported function in
$\operatorname{Lip}(\varepsilon,\C)$. It is then not difficult to
see that $\Phi$ is indeed bilipschitz. The main result of this paper
identifies a class of non-smooth functions $\mu$ which determine
bilipschitz quasiconformal mappings $\Phi$\,.

\begin{teor} \label{T}
Let $\{\Omega_j\}$\,, $1 \le j \le N$\,, be a finite family of
disjoint bounded domains of the plane with boundary of
class~$C^{1+\ep}$, $0<\ep <1$, and let $\mu = \sum_{j=1}^N
\mu_j\,\chi_{\Omega_j},$ where $\mu_j$ is of class
$\operatorname{Lip}(\ep,\Omega_j).$ \, Assume in addition that
$\|\mu\|_\infty < 1$. Then the associated quasiconformal mapping
$\Phi$ is bilipschitz.
\end{teor}

  Notice that the boundaries of the $\Omega_j$ may touch, even on a set of positive
  length and, of course, $\mu$ may have jumps on the boundary of some $
\Omega_j.\,$ In particular, if we only have one domain and $\mu$ is
constant we obtain the following corollary.

\begin{co} \label{C1}
If $\Omega$ is a bounded domain of the plane with boundary of
class $C^{1+\ep}$, $0<\ep <1$, and $\mu= \lambda \,\chi_\Omega$,
where $\lambda$ is a complex number such that $|\lambda|< 1$, then
the associated quasiconformal mapping~$\Phi$ is bilipschitz.
\end{co}

If $\Omega$ is a disc then Corollary 1 reduces to the fact that
$\Phi$ can be computed explicitly and that one can check by direct
inspection that is bilipschitz. If $\Omega$ is a square $Q$, then
one can show that the mapping $\Phi$ associated to
$\lambda\,\chi_Q$ is not Lipschitz for some $\lambda$ of modulus
less than $1$, so that the Corollary and thus the Theorem are
sharp as far as the smoothness of the boundaries of the $\Omega_j$
is concerned.

Recall that a $\mu$-quasi-regular function on a domain $D$ is a
complex function $f$ in $W^{1,2}_{\text{loc}}(D)$ satisfying (1),
with $\Phi$ replaced by $f$, almost everywhere in $D$. By Stoilow's
factorization theorem,\, $f = h \circ \Phi$ \,\,for some holomorphic
function $h$ on $\Phi(D).$  From the Theorem we then conclude that
$f$ is locally Lipschitz on $D$\,. This improves on Mori's Theorem,
which asserts that, for general $\mu$\,, $f$ is locally in
 $\operatorname{Lip}\alpha$ for $\alpha =
 \frac{1-\|\mu\|_\infty}{1+\|\mu\|_\infty}<1$\,. Thus, from the perspective of PDE, the Theorem
 may also be viewed as a regularity result for the Beltrami
 equation.

 The Beltrami equation is intimately
 related to second order elliptic equations in divergence form ot the type
 \begin{equation}\label{eq1bis}
\operatorname{div}\,(A \,\nabla u) =0\,,
 \end{equation}
 where $A=A(z)$ is a $2\times 2$  symmetric elliptic matrix with bounded
 measurable coefficients and determinant $1$ \,(see \cite[Chapter 13]{AIM}). Indeed, the real
 and imaginary parts of a solution to the Beltrami equation satisfy
  \eqref{eq1bis}, where the entries of the matrix $A$ are given explicitly in terms of the Beltrami coeficient.
  Conversely, given a solution $u$
 of \eqref{eq1bis}, one may find a solution of an appropiate
 Beltrami equation whose real part is $u.$
Thus for regularity issues one
 can work indistinctly with the Beltrami equation or with equation
 \eqref{eq1bis}. The proof of the Theorem gives, in particular,
  the following regularity result for solutions of equation \eqref{eq1bis}.

\begin{co}\label{C2}
Let $\Omega_j$\,, $1 \le j \le N$\,, be a finite family of
disjoint bounded domains of the plane with boundary of
class~$C^{1+\ep}$, $0<\ep <1$, and assume that each $\Omega_j$ is
contained in bounded domain $D$ with boundary of
class~$C^{1+\ep}$. Let $A=A(z),\,\, z \in D,$\, a $2\times 2$
symmetric elliptic matrix with determinant $1$ and entries
supported in $ \cup_{j=1}^N \Omega_j$ and belonging to
$\operatorname{Lip}(\ep,\Omega_j),$\,\, $1 \le j \le N .$ Let $u$
be a solution of equation \eqref{eq1bis} in $D.$ Let $D_\delta$
stand for the set of points in $D$ at distance greater than
$\delta$ from the boundary of $D.$ Then $\nabla u \in
\operatorname{Lip}(\varepsilon',\Omega_j \cap D_\delta),$\,\, for
$0 <\varepsilon' < \varepsilon ,$  and $1 \le j \le N .$ In
particular, $\nabla u \in L^\infty(D_\delta)$ and $u$ is a locally
Lipschitz function in $D.$
\end{co}

The main point of the corollary above is that each solution of
\eqref{eq1bis} is locally Lipschitz in $D,$ \, while the classical
De Giorgi-Nash Theorem gives only that $u$ satisfies locally a
Lipschitz condition of order $\alpha,$ \, for some $\alpha$
satisfying $0 <\alpha <1.$ See section 8 for an extension to more
general domains, which may have cusps.

There is some overlapping here with previous results by Li and
 Vogelius (\cite{LiV}) and Li and Nirenberg (\cite{LiN}). See at the end of the introduction for more
 about that.

 Another application of our Theorem concerns removability problems.
 There has recently been a renewed interest in gaining a better understanding
 of the nature of  removable sets for bounded quasi-regular functions (see
 \cite{ACMOU}, \cite{CFMOZ} and  \cite{CT}). Since bilipschitz mappings preserve removable
 sets for bounded holomorphic functions (\cite{T}), the Theorem
 immediately says that the removable sets for bounded $\mu$-quasi-regular
 functions, with $\mu$ as in the Theorem, are exactly the removable
 sets for bounded holomorphic functions.

If $\Omega$ is a domain, the $\operatorname{Lip} \ep$ norm of a
function $f$ on $\Omega$ is
\begin{equation}\label{eq2}
 \|f\|_{\ep} = \|f\|_{\ep,\Omega}  = \|f \|_{L^\infty(\Omega)}
+ \sigma_{\ep}(f)\,,
\end{equation}
where
\begin{equation}\label{eq2bis}
\sigma_{\ep}(f)= \sup \left\{ \frac{|f(z)-f(w)|}{|z-w|^{\ep}} :
z,w \in \Omega ,\, z \neq w \right\}\,.
\end{equation}

The main difficulty in proving the Theorem lies in the fact that no
smallness assumption is made on  $\sup_{1 \le j \le N}
\|\mu_j\|_{\ep,\Omega_j}$\,. In the same vein, Corollary 1 is much
more difficult to prove if $|\lambda|$ is close to~$1$. If one
assumes that $\|\mu_j\|_{\ep,\Omega_j}$ is small enough (depending
on $\Omega_j$) for each $j$, then the Theorem becomes easier.
Similarly, the Corollary 1 becomes easier under the assumption that
$|\lambda| \leq \ep_0(\Omega) \ll 1$\,. See a sketch of the argument
at the end of Section 2.

The scheme for the proof of the Theorem is inspired by a clever idea
of Iwaniec~\cite[p.~42--43]{I1} in the context of $L^p$ spaces,
which has been further exploited in \cite{AIS}. This idea brings
into play the index theory of Fredholm operators on Banach spaces
and, thus, compact operators. Our underlying Banach space is
$\operatorname{Lip}(\ep,\Omega)\,,\, \Omega$ a domain with boundary
of class~$C^{1+\ep}$, and on this space we estimate the Beurling
transform and its powers. We also show that the commutator between
the Beurling transform and certain functions
 is compact on appropriate larger Lipschitz spaces.

The Beurling transform is the principal value convolution operator

$$
Bf(z)= - \frac{1}{\pi}\, PV \int f(z-w)\frac{1}{w^2}\, dA(w)\,.
$$
The Fourier multiplier of $B$ is $\frac{\overline{\xi}}{\xi}$\,, or,
in other words,
$$
\widehat{Bf}(\xi) = \frac{\overline{\xi}}{\xi}\,\, \hat{f}(\xi)\,.
$$
Thus $B$ is an isometry on $L^2(\C)$\,.

 Our Main Lemma shows that
for each even smooth homogeneous Calder\'{o}n-Zygmund operator $T$ the
mapping
$$
T_\Omega(f)(z): = Tf(z)\,\,\chi_\Omega(z)\,,
$$
sends continuously $\Lipo$ into itself, where $\Omega$ is a
bounded domain with boundary of class $C^{1+\ep}\,.$  Throughout
the paper we understand that, for $f \in
\operatorname{Lip}(\varepsilon,\Omega),$ $Tf =
T(f\,\chi_\Omega)\,.$ The above boundedness result fails if $T$ is
not even. As a simple example, one may take as $T$ the Hilbert
transform and as $\Omega$ the interval $(-1,1)$.


The even character of $T$ is used in the proof of the Main Lemma in
the form
$$
T(\chi_D)\,\chi_D = 0\, ,\quad \text{for each disc}\,\, D\,,
$$
which should be understood as a local version of the global
cancellation property~$T(1)=0$ common to all smooth homogeneous
Calder\'{o}n-Zygmund operators. This was proved by Iwaniec for the
Beurling transform in \cite{I2}\,.

In Section 2 we present a detailed sketch of the proof and we
introduce the lemmas required. In Section 3 the Main Lemma is
proved in $\Rn$\,. Section 4 deals with commutators. We compare in
Section 5 the operator $B_\Omega^n$ with $ \chi_\Omega\,B^n$\,,\,
$\Omega$ a bounded domain with boundary of class~$C^{1+\ep}$\,,
and we show that the difference is compact on $\Lip$, for
$0<\ep'<\ep$. In Section 6 one completes the proof of the Theorem
for the case of one domain. Section 7 contains the reduction to
the one domain case. In section 8 we present an extension of the
Theorem to what seems to be its more natural setting, namely
domains with cusps whose boundary is of class $C^{1+\varepsilon}$
off the set of cusps. Applications to the regularity theory of
solutions of equation \eqref{eq1bis} is this setting are also
mentioned.

After a first version of the paper was completed,  Daniel Faraco
brought to our attention the work of Li and Vogelius \cite{LiV} (and
\cite{LiN}), which deals with Lipschitz regularity for the equation
\eqref{eq1bis} in $\Rn$ for matrices with entries satisfying a
Lipschitz condition of order $\varepsilon$ on finitely many disjoint
domains with boundary of class $C^{1+\varepsilon}$, but with
possible jumps across the boundaries. The closures of the domains
were disjoint and a main point was to obtain gradient estimates
independent of the mutual distances between the closed domains. This
is not an issue for our methods, which even allow touching domains.
Moreover, as stated before in Corollary 2, for each solution of
\eqref{eq1bis} we obtain a regularity of class
$C^{1+\varepsilon'}$\,, \,\, for each $\varepsilon' <
\varepsilon,$\, in each domain. This is almost the expected best
possible result, namely $C^{1+\varepsilon}$ in each domain. In
\cite{LiV} there is a more substantial loss, due to the techniques
employed. On the other hand, the setting in \cite{LiV} and
\cite{LiN} is more general, in the sense that one works in $\Rn$,
there is no restriction on the determinant of the matrix  and also
non-homogeneous terms are considered.

We also learnt from Antonio C\'{o}rdoba that the regularity theory of
the Euler equation in 2D, in particular the regularity theory of
vortex patches, makes broad use of even Calder\'{o}n-Zygmund operators,
the Beurling transform in particular. We then became aware of the
article \cite{D}, in which one also proves the Main Lemma. However,
the proof there is different, and certainly not as much in the
Calder\'{o}n-Zygmund tradition as ours.

\section{Sketch of the proof}

First of all, there is a standard factorization method in
quasiconformal mapping theory that reduces the Theorem to the case
of only one domain $\Omega \, (N=1)\,.$  The argument is presented
in detail in section 7. Then, from now on we will assume that $\mu$
vanishes off some domain $\Omega$ with boundary of class~$C^{1+\ep}$
 and that $\mu \in \operatorname{Lip}(\ep,\Omega)\,.$

 As is well known, $\Phi$ is given explicitly by the formula \cite
{Ah}
$$
\Phi(z) = z + C(h)(z)\,,
$$
where
$$Ch(z) = \frac{1}{\pi}\, \int h(z-w)\frac{1}{w}\, dA(w)$$
is the Cauchy transform of $h$\,. Recall the important relation
between the Cauchy and the Beurling transforms: $\partial C = B$,
$\partial= \frac{\partial}{\partial z}$. The function $h=
\overline{\partial}\,{\Phi}$ is determined by the equation
$$
(I-\mu\,B)(h) = \mu\,.
$$
As soon as we can invert the operator $I-\mu\,B$ on $\Lip$, for some
$\ep'$ satisfying $0<\ep'<\ep$\,,  then
$$
h= (I-\mu\,B)^{-1}(\mu)\,,
$$
and thus $h$ is in $\Lip$ and, in particular, is bounded on
$\Omega$. By the Beltrami equation \eqref{eq1}\,, $h$ vanishes on
$\C \setminus \overline{\Omega}$ and therefore $h$ is in
$L^\infty(\C)$\,. On the other hand,  $\partial\, \Phi = 1+ B(h)$.
By the Main Lemma\,, $B(h)$ is in  $L^\infty(\C)$ (see \eqref{eq8}
below)\,, and so $\Phi$~is a Lipschitz function on the plane.

Showing that $\Phi$ is bilipschitz still requires an argument.
Indeed, we have shown up to now that $\Phi$ is of class
$C^{1+\ep'}(\Omega)$ and thus its Jacobian is non-zero at each
point of $\Omega$ (\cite[Theorem 7.1, p.~233]{LV})\,. On the other
hand, $\Phi$ is conformal on $\C \setminus \overline{\Omega}$ and
thus the Jacobian is also non-zero there. However we cannot infer
immediately that the Jacobian is bounded below away from zero
either on $\Omega$ or on $\C \setminus \overline{\Omega}$. This is
proved in Section~\ref{sec6} and hence $\phi$~is bilipschitz.

It remains to prove that $I-\mu B$ is invertible on $\Lip$ for
each $\ep'$ with $0<\ep'<\ep$.  For $f$ in $\Lipo$ set
$$
B_\Omega(f)(z) = B(f)(z)\,\,\chi_\Omega(z)\,,
$$
where, as we said in the introduction, $B(f)$ stands for
$B(f\,\chi_\Omega).$  Following \cite[p.~48]{AIS} we define
$$
P_m = I+\mu B_\Omega+ (\mu B_\Omega)^2 + \cdots + (\mu
B_\Omega)^{m}\,,
$$
so that we have
\begin{equation}\label{eq3}
 (I-\mu B_\Omega) P_{n-1}=P_{n-1} (I-\mu B_\Omega)= I- (\mu
B_\Omega)^n = I- \mu^n B_\Omega^n +R \,,
\end{equation}
where $R= \mu^n B_\Omega^n - (\mu B_\Omega)^n$ can be easily seen to
be a finite sum of operators that contain as a factor the commutator
$K_0 =\mu B_\Omega- B_\Omega \mu$. Lemma 3 in Section 4 asserts that
$K_0$ is compact on $\Lip$ for each $\ep'$ less than $\ep$\,, so
that $R$ is also compact on $\Lip$. One would like to have now that
the operator norm of $\mu^n B_\Omega^n$ on $\Lip$ is small if $n$ is
large. Would this be so, then
 $I-\mu B_\Omega$ would be a Fredholm operator on $\Lip$. But it looks like a difficult task to
 obtain estimates
 for the operator norm of $B_\Omega^n$ better than the obvious exponential upper bound
 $\|B_\Omega\|^n$\,. We overcome this difficulty by finding an
 expression of the form
 \begin{equation}\label{eq4}
 B_\Omega^n(f) = B^n(f)\,\chi_\Omega+ K_n(f)\,,
 \end{equation}
where $K_n$ is compact on $\Lip$. This is done in Theorem 1 in
Section 5\,. Incidentally, in turns out  that $K_n =0$ when $\Omega$
is a disc, so that in this case $ B_\Omega^n(f)$ is exactly
$B^n(f)\,\chi_\Omega$ for each $n$\,.

Then \eqref{eq3} can be rewritten as
\begin{equation}\label{eq5}
(I-\mu B_\Omega) P_{n-1}=P_{n-1} (I-\mu B_\Omega)= I-\mu^n B^n +
S\,,
\end{equation}
where $S$ is compact on $\Lip$.

The kernel of $B^n$ may be computed explicitly, for instance  via a
Fourier transform argument \cite [p.~73]{St}, and one obtains
$$
b_n(z)= \frac{(-1)^n n}{\pi}\,\frac{\bar z^{n-1}}{z^{n+1}}\,.
$$
Thus the Calder\'{o}n-Zygmund constant of $b_n$\,, namely,
$$
\|b_n(z)\,|z|^2\|_\infty + \| \nabla b_n(z)\,|z|^3\|_\infty \,,
$$
is less than $C n^2$, where $C$ is a positive constant. Hence, by
the Main Lemma
$$
\| \mu^n B^n(f)\|_{\ep',\Omega} \leq C\, n^3 \,\|\mu\|_\infty^n
\,\|\mu\|_{\ep',\Omega}\, \|f\|_{\ep',\Omega}\,,
$$
which tells us that the operator norm of $\mu^n B^n$ as an operator
on $\Lip$ is small for large $n$\,. Therefore $I-\mu B_\Omega$ is a
Fredholm operator on $\Lip$.

 Clearly $I- t \,\mu B_\Omega$\,,\, $0 \leq t \leq 1$ is a continuous path from
the identity to $I-\mu B_\Omega$\,. By the index theory of Fredholm
operators on Banach spaces (e.g.~\cite{Sch}), the index is a
continuous function of the operator. Hence $I-\mu B_\Omega$ has
index $0$. On the other hand $I-\mu B_\Omega$ is injective, because
if $f=\mu B_\Omega(f)$, then $\|f\|_2 \leq \|\mu\|_\infty\,
\|B_\Omega(f)\|_2 \leq \|\mu\|_\infty\, \|B(f)\|_2 =
\|\mu\|_\infty\, \|f\|_2$\,, which is possible only if $f=0$\,. Thus
$I-\mu B_\Omega$ is invertible on $\Lip$.

As we mentioned before, the proof of the Theorem simplifies if
$\|\mu\|_{\ep, \Omega} $ is assumed to be less than a small number
$\delta_0 = \delta_0(\Omega)$\,. In this case one can invert $I-\mu
B$ by a Neumann series and get $h= \sum_{n=0}^\infty (\mu
B)^n(\mu)\,.$ By the Main Lemma $B_\Omega$ is bounded on $\Lipo$\,.
Denote by $\|B_\Omega\|$ its operator norm and assume that
$\|\mu\|_{\ep,\Omega} < (2\,\|B_\Omega\|)^{-1}$. Then
$$
\|h\|_{\ep,\Omega} \leq \sum_{n=0}^\infty \|\mu
\|^{n+1}_{\ep,\Omega} \, \|B_\Omega\|^n \le 2\,\|\mu\|_{\ep,\Omega}
< \|B_\Omega\|^{-1}\,.
$$
But is also part of the Main Lemma that
$$
\|B(h)\|_{L^\infty(\C)} \le C(\Omega)\, \|h\|_{\ep,\Omega}\,.
$$
Hence, if we also assume that $2\,\|\mu\|_{\ep,\Omega}\,\, C(\Omega)
< 1$, we have $\|B(h)\|_{L^\infty(\C)} < 1\,.$ Thus
$$
\|\overline{\partial}\,\Phi\|_{L^\infty(\C)} <
\|B_\Omega\|^{-1}\quad \text{and}\quad \|\partial
\,\Phi\|_{L^\infty(\C)} \leq 1+ \|B(h)\|_{L^\infty(\C)} \leq 2\,,
$$
and so $\Phi$ is a lipschitz function\,. That $\Phi$ is bilipschitz
follows from
$$
|\partial \,\Phi(z)|= |1+ B(h)(z)| \geq 1- \|B(h)\|_{L^\infty(\C)}
> 0 \,, \quad z \in \C \setminus \fr \,.
$$

\section{The Main Lemma}
In this section we move to $\Rn$. We say that a bounded domain
$\Omega \subset \Rn$ has a boundary of class $C^{1+\ep}$ if $\fr$ is
a $C^1$ hyper-surface whose unit normal vector satisfies a Lipschitz
condition of order $\ep$ as a function on the surface. To state an
alternative condition, for $x = (x_1, \dots ,x_n)\in \Rn $ we use
the notation $x = (x',x_n)$\,, where $x'= (x_1,\dots,x_{n-1})$\,.
Then $\Omega$ has a boundary of class $C^{1+\ep}$ if for each point
$a \in
\partial\,\Omega$ one may find a ball $B(a,r)$ and a function $x_n =
\varphi(x')$, of class $C^{1+\ep}$, such that, after a rotation if
necessary, $\Omega \cap B(a,r)$ is the part of $B(a,r)$ lying below
the graph of $\varphi$\,. Thus we get
\begin{equation}\label{eq6}
\Omega \cap B(a,r)= \{x \in B(a,r) : x_n <
\varphi(x_1,\dots,x_{n-1})\}\,.
\end{equation}

A smooth (of class $C^1$) homogeneous Calder\'{o}n-Zygmund operator is a
principal value convolution operator of type
\begin{equation}\label{eq6bis}
T(f)(x)= PV \int f(x-y)\,K(y) \,dy \,,
\end{equation}
where
$$
K(x)= \frac{\omega(x)}{|x|^{n}}\,,\quad x \neq 0\,,
$$
$\omega(x)$ being a homogeneous function of degree~$0$\,,
continuously differentiable on $\Rn \setminus\{0\}$ and with zero
integral on the unit sphere\,. The maximal singular integral
associated to $T$ is
$$
T^{\star}f(x)= \sup_{\delta > 0} | T^{\delta}f(x)|, \quad x \in
\Rn\,,
$$
where
$$
T^{\delta}f(x)= \int_{| y-x| > \delta} f(x-y) K(y) \,dy\,.
$$
The Calder\'{o}n-Zygmund constant of the kernel of $T$ is defined as
$$
\| T \|_{CZ} = \|K(x)\,|x|^n \|_\infty + \|\nabla
K(x)\,|x|^{n+1}\|_\infty\,.
$$
The operator $T$ is said to be even if the kernel is even, namely,
if $\omega(-x)=\omega(x)\,,$ \,for all $ x \neq 0\,.$

We are now ready to state our main lemma. The definition of the norm
in $\Lipo$ is as in \eqref{eq2}. As we explained in the previous
section, we need the precise form of the constant in the inequality
below.

\begin{ML}
Let $\Omega$ be a bounded domain with boundary of class $C^{1+\ep}$,
$0<\ep<1$, and let $T$ be an even smooth homogeneous
Calder\'{o}n-Zygmund operator. Then $T$ maps $\Lipo$ into $\Lipo$, and
 $T$ also maps $\Lipo$ into $\operatorname{Lip}(\ep,\Omega^c)$. In fact, one
has the inequalities
$$
\| Tf \|_{\ep,\Omega} \leq C\, \|T\|_{CZ}\,\|f\|_{\ep,\Omega}\,,
$$
and
$$
\| Tf \|_{\ep,\Omega^c} \leq C\, \|T\|_{CZ}\,\|f\|_{\ep,\Omega}\,,
$$
 where $C$ is a constant depending only on $n$\,, $\ep$
and $\Omega$\,.
\end{ML}

\begin{proof} We choose a positive $r_0 = r_0(\Omega)$ small enough so that a series of properties
that will be needed along the proof are satisfied. The first one is
that for each $a \in \partial\,\Omega$\,, which we can assume to be
$a=0$\,, we have \eqref{eq6}. After a rotation we may assume that
the tangent hyperplane to $\partial\,\Omega$ at $0$ is $x_n=0$\,. We
take $r_0$ so small that
\begin{equation}\label{eq7}
|\varphi(x')| \leq C |x'|^{1+\ep}\,, \quad |x'| < r_0\,,
\end{equation}
for some positive constant $C$ depending only on $\Omega$\,. We
claim that
\begin{equation}\label{eq8}
T^*f (x) \leq C \, \|T\|_{CZ} \,\|f\|_\ep\,,\quad x \in \Rn\,.
\end{equation}
The proof of \eqref{eq8} is a technical variation of the proof of
Lemma 5 in \cite{MOV}\,.  We have
\begin{equation*}
\begin{split}
T^\delta(f)(x)&= \int_{\delta < |y-x|< r_0} f(y)\, K(x-y)\,dy +
\int_{r_0 < |y-x|} \dotsi\\*[5pt] &= I_\delta + II \,.
\end{split}
\end{equation*}
Clearly,
$$
|II| \leq  \int_{r_0 < |y-x|} |f(y)|\, |K(x-y)|\,dy \leq
r_0^{-n}\,|\Omega|\,\|T \|_{CZ}\, \|f\|_\infty \,.
$$
To deal with the term $I_\delta$ we write
\begin{equation*}
\begin{split}
I_\delta&= \int_{\delta < |y-x|< r_0} \chi_\Omega(y)\,(f(y)-f(x))\,
K(x-y)\,dy\\*[5pt] &\quad+ f(x)\int_{\delta < |y-x|< r_0}
\chi_\Omega(y)\, K(x-y)\,dy
\\*[5pt]
& = III_\delta+ f(x)\,IV_\delta \,,
\end{split}
\end{equation*}
and we remark that $III_\delta$ can easily be estimated as follows
\begin{equation*}
\begin{split}
|III_\delta| & \leq \|f\|_{\ep}\,\int_\Omega |y-x|^\ep
|K(x-y)|\,dy \\
& \leq C\, \|f \|_{\ep}\, \|T\|_{CZ}\, \int_\Omega
|y-x|^{-n+\ep}\,dy \\ & \leq C(\ep)\,(\text{diam} \Omega)^\ep \,\|f
\|_{\ep}\, \|T\|_{CZ}\,.
\end{split}
\end{equation*}
Taking care of $IV_\delta$ is not so easy. Assume first that $x =0$
is in $\fr$\,. Without loss of generality we may also assume  that
the tangent hyperplane to $\fr$ at $0$ is $\{x_n =0\}\,$ (see Figure
1)\,.

\begin{figure}[ht]
\begin{center}
\includegraphics{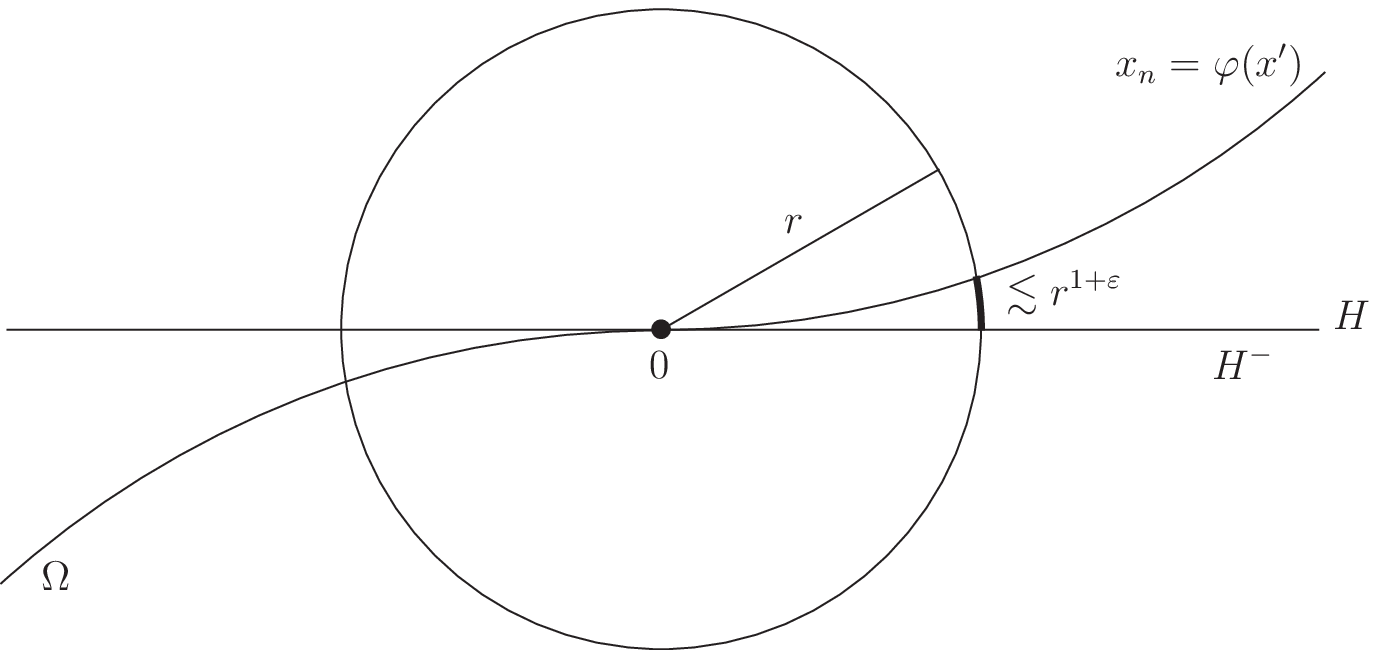}\\*[7pt]
Figure 1
\end{center}
\end{figure}


\noindent Let $H_-$ be the half space $\{x_n < 0 \}$\,. Take
spherical coordinates $y=r\,\xi$ with $0 \leq r$ and $|\xi|=1$. Then
\begin{equation}\label{eq9}
IV_\delta  =  \int_\delta^{r_0} \left( \int_{A(r)}
\omega(\xi)\,d\sigma(\xi)\right)\frac{dr}{r} \,,
\end{equation}
where
$$
A(r)= \{ \xi : |\xi| =1 \text{ and } r \xi \in \Omega \}\,,
$$
and $\sigma$ is the surface measure on the unit sphere $U$\,. Since
$K$ is even,
$$
0= \int_U \omega(\xi)\,d\sigma(\xi) = 2\, \int_{U \cap H_-}
\omega(\xi)\,d\sigma(\xi)\,.
$$
Thus
$$
\int_{A(r)} \omega(\xi)\,d\sigma(\xi) =  \int_{A(r)\setminus (U \cap
H_-)} \omega(\xi)\,d\sigma(\xi) - \int_{(U \cap H_-)\setminus A(r)}
\omega(\xi)\,d\sigma(\xi)\,,
$$
and so
$$
\left|\int_{A(r)} \omega(\xi)\,d\sigma(\xi)\right| \leq C\,\|
T\|_{CZ} \, \left(\sigma(A(r) \setminus (U \cap H_-) ) + \sigma((U
\cap H_-)\setminus A(r))\right)\,.
$$
By \eqref{eq7}, we obtain
\begin{equation*}
\sigma(A(r) \setminus (U \cap H_-) ) + \sigma((U \cap H_-)\setminus
A(r)) \leq C\,r^\ep \,,
\end{equation*}
which yields, by \eqref{eq9},
\begin{equation*}
|IV_\delta| \leq C\,\|T\|_{CZ}\,.
\end{equation*}

Take now $x \in \Rn \setminus \fr$\,.  Denote by $\delta_0 $ the
distance from $x$ to $\fr$ and let $x_0$ be a point in $\fr$ where
such distance is attained. Set
$$
A= \{y \in \Omega : \delta_0 < |y-x| < r_0 \}
$$
and
$$
A_0= \{y \in \Omega : \delta_0 < |y-x_0| < r_0 \}\,.
$$
We compare $IV_\delta$ to the expression we get replacing $x$ by
$x_0$ and $\delta$ by $\delta_0$ in the definition of $IV_\delta$.
For $\delta \leq \delta_0$ we have, by the standard cancellation
property of the kernel,
$$
\int_{\delta < |y-x|< r_0} \chi_\Omega(y)\, K(x-y)\,dy =
\int_{\delta_0 < |y-x|< r_0} \chi_\Omega(y)\, K(x-y)\,dy
$$
and then
\begin{equation*}
\begin{split}
\biggl|\int_{\delta < |y-x|< r_0} &\chi_\Omega(y)\, K(x-y)\,dy -
\int_{\delta_0 < |y-x_0|< r_0} \chi_\Omega(y)\, K(x_0-y)\,dy \biggr|
\\*[7pt] &= \left|\int_{A}  K(x-y)\,dy - \int_{A_0} K(x_0-y)\,dy
\right|  \\*[7pt] &\leq \int_{A\cap A_0}
|K(x-y)-K(x_{0}-y)|\,dy\\*[7pt] &\quad+ \left|\int_{A \setminus A_0}
\!\chi_\Omega(y)\, K(x-y)\,dy \right| +
 \left| \int_{A_0 \setminus A} \chi_\Omega(y)\, K(x_0-y)\,dy  \right|  \\*[7pt]
&= J_1+J_2+J_3\,.
\end{split}
\end{equation*}
If $y \in A \cap A_0$, then
$$
|K(x-y)-K(x_0-y)| \leq C \,\|T\|_{CZ}\,
\frac{|x-x_0|}{|y-x|^{n+1}}\,.
$$
Hence
$$
J_1 \leq C \,\|T\|_{CZ}\, |x-x_0|\, \int_{|y-x|> \delta_0}
\frac{dy}{|y-x|^{n+1}} \leq C\,\|T\|_{CZ}\,.
$$
To estimate $J_2$ observe that
$$
A \setminus A_0 = \left(A \cap B(x_0,\delta_0)\right) \cup \left(A
\cap (\Rn \setminus B(x_0,r_0))\right)\,.
$$
Assume for the moment that $\delta_0 \leq r_0/2$\,.  Now, it is
obvious that if $|y-x_0|\geq r_0$, then $|y-x| \geq r_0/2$, and so
$$
J_2 \leq \|T \|_{CZ} \left( \int_{|y-x_0|< \delta_0}
\frac{dy}{\delta_0^n} + \int_\Omega \frac{2^n}{r_0^n} \,dy \right)
\leq C\,\|T \|_{CZ} \,.
$$
A similar argument does the job for $J_3$.

If $\delta_0 \geq r_0/2$\,, then  the estimate of $T^\delta (f)(x)$
is straightforward. Indeed, by the cancellation of the kernel we may
assume that $\delta \geq \delta_0$ and so
$$
|T^\delta (f)(x)| \leq \|f \|_\infty \,\|T\|_{CZ} \int_{|y-x|>
\delta} \chi_\Omega(y)\,\frac{1}{|y-x|^n}\, dy \leq
\frac{2^n}{r_0^n}\,|\Omega|\, \|f \|_\infty \,\|T\|_{CZ}\,.
$$
This completes the proof of \eqref{eq8}\,.

\pagebreak

Our next task is to estimate the semi-norm $\sigma_\ep(T(f))$ on
$\Omega$\,. For this we need a lemma, which should be viewed as a
manifestation of the extra cancellation enjoyed by even kernels.
Notice that no smoothness assumptions on the kernel are required.
The lemma is known for the Beurling transform (\cite[p.~389]{I2}).

\begin{lemma}\label{cancel}
Let $K(x)= \frac{\omega(x)}{|x|^{n}}$, where $\omega$ is an even
homogeneous function of degree~$0$, integrable on the unit sphere
and with vanishing integral there. Let $T$ be the associated
Calder\'{o}n-Zygmund operator defined by \eqref{eq6bis}. Then
$$
T(\chi_B)\,\chi_B = 0\,,\quad \text{for each ball}\,\, B\,.
$$
\end{lemma}
\begin{proof}
Assume, without loss of generality that $B$ is the unit ball. Fix a
point $x$ in~$B$\,. Then
\begin{equation*}
\begin{split}
T(\chi_B)(x)  & = \lim_{\ep \rightarrow 0} \int_{B \cap
B^c(x,\ep)} K(x-y)\,dy  \\
 & =\int_{B \cap B^c(x,1-|x|)}
K(x-y)\,dy = \int_{B^c \cap B(x,1+|x|)} K(x-y)\,dy\,.
\end{split}
\end{equation*}
Expressing the latest integral above in polar coordinates $y=x+r\xi$
centered at $x$,\, we get
$$
T(\chi_B)(x) = \int_{|\xi|=1} \int^{1+|x|}_{r(x,\xi)} \frac{dr}{r}\,
\omega(\xi)\,d\sigma(\xi) = \int_{|\xi|=1}  \log
\left(\frac{1+|x|}{r(x,\xi)}\right)\, \omega(\xi)\,d\sigma(\xi)\, .
$$
The lower value $r(x,\xi)$ is determined as shown in Figure 2.
\begin{figure}[ht]
\begin{center}
\includegraphics{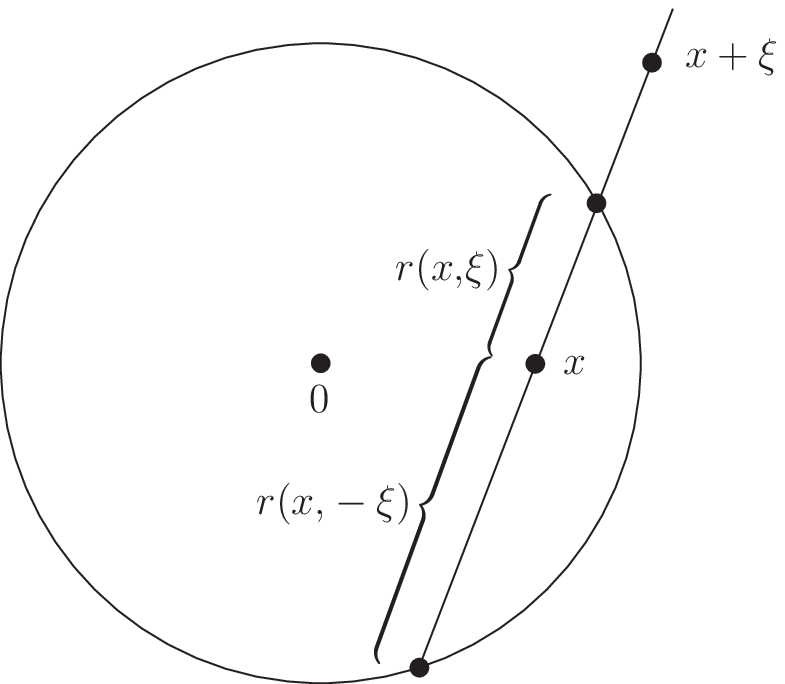}\\*[7pt]
Figure 2
\end{center}
\end{figure}

\noindent Since $\omega$ has zero integral on the unit sphere,
$$
T(\chi_B)(x) = \int_{|\xi|=1}  \log
\left(\frac{1}{r(x,\xi)}\right)\, \omega(\xi)\,d\sigma(\xi)\, .
$$
Set $U^+$ be the half of the unit sphere above the hyperplane
$\{x_n=0\}\,$. Since $\omega$ is even,
$$
T(\chi_B)(x) = \int_{U^+}  \log
\left(\frac{1}{r(x,\xi)\,r(x,-\xi)}\right)\,
\omega(\xi)\,d\sigma(\xi)\,.
$$
Now, the points $0$\,,\,$x$, $x+\xi$ and $x-\xi$ lie in a plane that
intersects the unit sphere in a circumference. It is clear from
Figure 2 that the product $r(x,\xi)\,r(x,-\xi)$ is the power of $x$
with respect to that circumference and thus it does not depend on
$\xi$ (in fact it is exactly $1-|x|^2)$ \,. Since the integral of
$\omega$ on each semi-sphere is zero the proof is complete.
\end{proof}

We are now ready to estimate the semi-norm $\sigma_\ep(T(f))$\,. We
deal first with the case $f=1$. The general case follows from this
by the $T(1)$-Theorem for Lipschitz spaces on spaces of homogeneous
type \cite {W} (see also \cite{Ga} and \cite{GG} for the
non-doubling case)\,. The conditions on the kernel required in \cite
{W} (and in \cite {Ga}, \cite {GG}) are implied by the fact that
$T^*(\chi_\Omega) \in L^\infty(\Omega) $\,, which we proved
before\,. However, in our particular setting the reduction to $f=1$
is elementary and will be discussed afterwards for the sake of
completeness\,.

We want to prove that
\begin{equation}\label{eq10}
|T(\chi_\Omega)(x)-T(\chi_\Omega)(y) | \leq C\, |x-y|^\ep\,,\quad
x\,,y \in \Omega\,.
\end{equation}
Fix $x$ and $y$ in $\Omega$\,. Changing notation if necessary,  we
can assume that $\operatorname{dist}(x,\fr) \leq
\operatorname{dist}(y,\fr)$\,. We may also assume, without loss of
generality, that $\operatorname{dist}(x,\fr) \leq r_0 / 4$\,.
Otherwise we have
$$
|T(\chi_\Omega)(x)-T(\chi_\Omega)(y) | \leq C \, |x-y|\,\|\nabla
T(\chi_\Omega)\|_{L^\infty(\Omega_{0})}\,,
$$
where $\Omega_0 = \{z \in \Omega : \operatorname{dist}(x,\fr) \geq
\R / 4 \}$\, and $C$ depends only on $\Omega$. Notice that
$T(\chi_\Omega) \in C^1(\Omega)$\,, because
$$T(\chi_\Omega)= T(1-\chi_{\C \setminus \Omega}) = - T(\chi_{\C
\setminus \Omega})$$ and the kernel of $T$ is continuously
differentiable off the origin. Indeed, for some constant depending
only on $n, \R $ and  $\Omega$\,, we have
$$ \|\nabla T(\chi_{\C
\setminus \Omega})\|_{L^\infty(\Omega_{0})} \leq C\,\|T\|_{CZ}\,.$$

We may also assume, without loss of generality, that $|x-y| \le \R/
4$\,, because, otherwise,
$$
|T(\chi_\Omega)(x)-T(\chi_\Omega)(y) | \leq
\frac{8}{\R}\,\|T(\chi_\Omega)\|_\infty\,|x-y|\,.
$$

Having settled these preliminaries we proceed to the core of the
proof of \eqref{eq10}\,.  We may assume that the point of $\fr$
nearest to $x$ is the origin. Let $B$ be the ball with center
$(0,\dots,0,-\R)$ and radius $\R$, so that $\partial\,B$ is tangent
to $\fr$ at $0$\,. Let $S$ stand for the set $(\Omega \setminus B)
\cup (B \setminus \Omega)$\,. The central idea in the proof of the
Main Lemma is to use the extra cancellation of even Calder\'{o}n-Zygmund
operators  via Lemma 1 to write
$$
|T(\chi_\Omega)(x)-T(\chi_\Omega)(y) | =
|T(\chi_S)(x)-T(\chi_S)(y)|\,.
$$

\pagebreak

\noindent The obvious advantage is that $S$ is a region which is
``tangential"  to $\fr$ at $0$\,, and hence small\,.  By \eqref{eq7}
we may take $\R$ so small that for some constant $C$ depending only
on $\Omega$
\begin{equation}\label{eq10bis}
|(z,\overrightarrow{n})| \leq C \,|z'|^{1+\ep}\,,\quad z \in S \cap
B(0,\R)\,,
\end{equation}
where $\overrightarrow{n}$ stands for the inward unit normal vector
to $\fr$ at $0$ and $(\; ,  \; ) $ denotes the scalar product  in
$\Rn$. Thus, if $\R$ is small enough,
\begin{equation*}
|(z,\overrightarrow{n})| < \frac{1}{\sqrt{2}}\,|z|\,,\quad z \in S
\cap B(0,\R)\,.
\end{equation*}

We distinguish two cases according to whether the position of $x$
and $y$ relative to $\fr$ is non-tangential or tangential\,. To make
this precise we introduce the cone~$\Gamma$ with vertex $0$ and
amplitude $\pi/2$\,, namely,
\begin{equation}\label{eq11}
\Gamma = \{z \in \C : (z,\overrightarrow{n})\geq
\frac{1}{\sqrt{2}}\,|z|\}\,.
\end{equation}
Clearly, if $\R$ is chosen small enough, then the part of the cone
near $0$ is contained in $\Omega$ and in $B$\,. More precisely,
$$ \Gamma \setminus \{0\} \cap B(0,\R) \subset \Omega \cap B
\cap B(0,\R)\,. $$  See Figure 3.

\vspace*{-.4cm}

\begin{figure}[h]
\begin{center}
\includegraphics{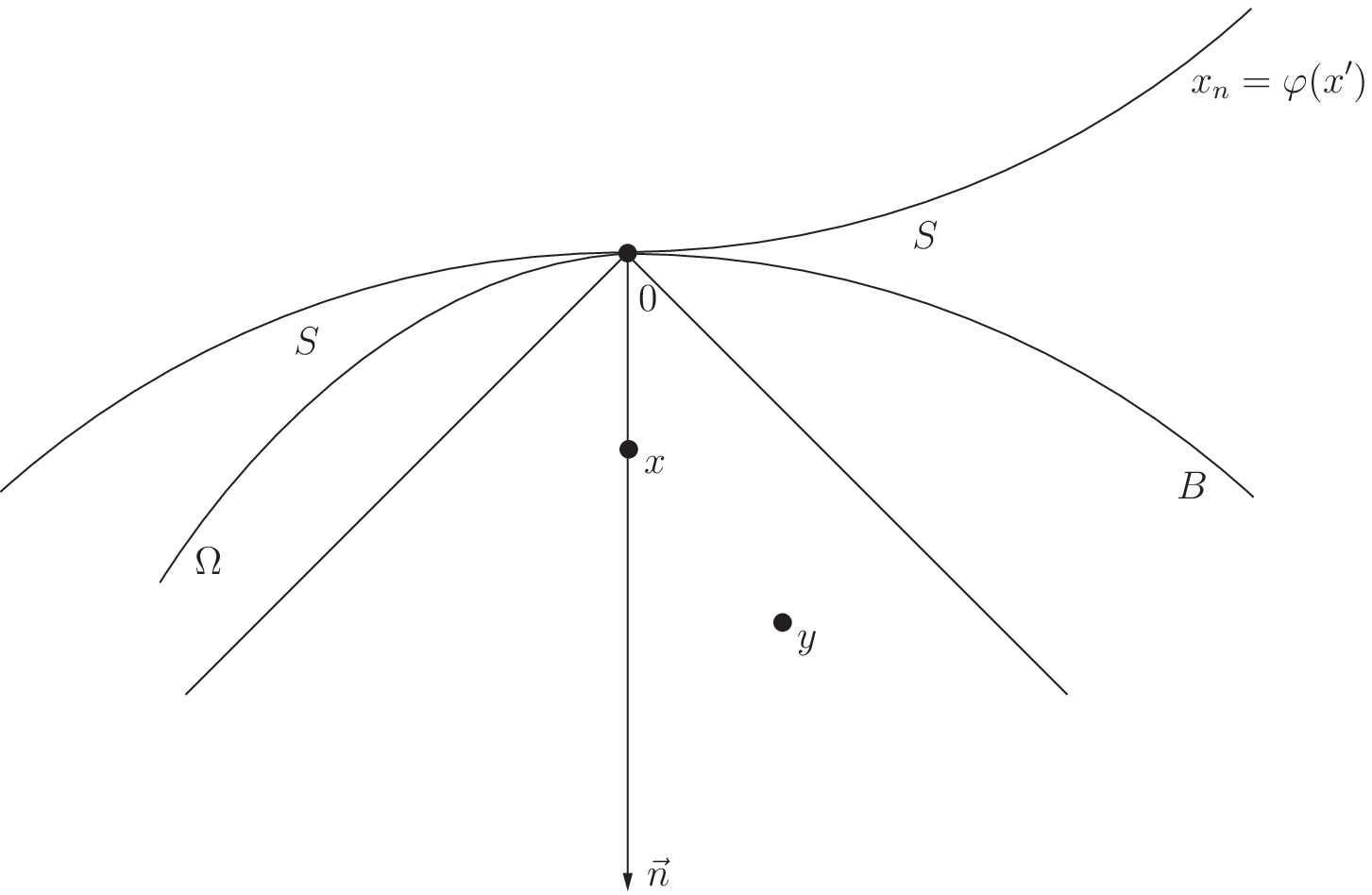}\\*[7pt]
Figure 3
\end{center}
\end{figure}

\noindent {\it Case 1:} $x$ and $y$ are in non-tangential position,
that is, $x$ and $y$ belong to $\Gamma$\,.  We have
\begin{equation*}
\begin{split}
|T(\chi_S)(x)-T(\chi_S)(y)| & \leq \left|\int_{S \cap B^c(x,\R)}
K(x-z)\,dz -  \int_{S \cap B^c(y,\R)} K(y-z)\,dz \right| \\*[7pt]
&\quad +  \left |\int_{S \cap B(x,\R)} K(x-z)\,dz - \int_{S \cap
B(y,\R)} K(y-z)\,dz \right| \\*[7pt] & = I+II\,.
\end{split}
\end{equation*}
Split $I$ into three terms as follows
\begin{equation*}
\begin{split}
I & \leq \left|\int_{S \cap B^c(x,\R)\cap B(y,\R)} K(x-z)\,dz \right
|+  \left|\int_{S \cap B^c(y,\R)\cap B(x,\R)} K(y-z)\,dz \right|
\\*[7pt] &\quad +  \left |\int_{S \cap B^c(x,\R)\cap B^c(y,\R)}
\left(K(x-z)-K(y-z)\right)\,dz \right| \\*[7pt] & = I_1+I_2 +I_3 \,.
\end{split}
\end{equation*}
The terms $I_1$ and $I_2$ are estimated in the same way. For
instance, for $I_1$, we get
\begin{equation*}
\begin{split}
|I_1| & \leq \int_{B^c(x,\R)\cap B(y,\R)} \frac{C}{|x-z|^n} \,dz
\leq \frac{C}{\R^n}\, |B^c(x,\R)\cap B(y,\R)|\\*[7pt] & \leq
\frac{C}{\R^n}\,|x-y|\,\R^{n-1} = \frac{C}{\R}\,|x-y|\,,
\end{split}
\end{equation*}
where in the latest inequality we used that $|x-y| \leq \R \,.$

The term $I_3$ is controlled by a gradient estimate, namely,
$$
I_3 \leq \int_{B^c(x,\R)} C\, \frac{|x-y|}{|x-z|^{n+1}}\,dz \leq
\frac{C}{\R}\,|x-y|\,.
$$

The more difficult term $II$ is not greater than
\begin{equation*}
\begin{split}
& \left |\int_{S \cap B(x,2\,|x-y|)} K(x-z)\,dz \right|\\*[9pt]
&\qquad\quad +\left| \int_{S \cap B(y,2\,|x-y)} K(y-z)\,dz
\right|\\*[9pt]
 &\qquad\quad +
\left|\int_{S \cap B(x,\R)\cap B^c(x,2\,|x-y|)} K(x-z)\,dz - \int_{S
\cap B(y,\R)\cap B^c(y,2\,|x-y|)} K(y-z)\,dz \right|
\\*[9pt]
&\qquad = II_1+ II_2+III\,.
\end{split}
\end{equation*}

\pagebreak

\noindent Estimating the three terms above requires a simple lemma.

\begin{lemma}\label{Lemma 2}
If $\R$ is small enough, then one has
$$
|w-z| \geq C\,|z|\,, \quad w \in \Gamma \cap B(0,\R)\,, \quad z \in
S \cap B(0,\R)\,,
$$
for  $C = \left(2 (1+\sqrt{2})\right)^{-1}$\,.
\end{lemma}
\begin{proof}
According to the definition of the cone $\Gamma$ and by \eqref{eq7}
\begin{equation*}
\begin{split}
|z| & \leq |z-w|+|w| \leq |z-w|+ \sqrt{2}\,
(w,\overrightarrow{n})\\*[5pt] & \leq (1+\sqrt{2})\,|z-w|+
\sqrt{2}\,|(z,\overrightarrow{n})|\\*[5pt] & \le
(1+\sqrt{2})\,|z-w|+ \sqrt{2}\,C\, |z|^{1+\ep }\\*[5pt] & \le
(1+\sqrt{2})\,|z-w|+ \sqrt{2}\,C\,\R^\ep \, |z|\,.
\end{split}
\end{equation*}
If $\R$ satisfies $\sqrt{2}\,C\,\R^\ep \le 1/2$\,, then $|z|\le
2(1+\sqrt{2})\,|z-w|$\,, which proves the lemma\,.
\end{proof}

To estimate the term $II_1$ we apply Lemma 2 to $w=x$ to obtain
$$
II_1 \le \int_{S \cap B(x,2\,|x-y|)} \frac{C}{|x-z|^n}\,dz \le C\,
\int_{S \cap B(x,2\,|x-y|)} \frac{dz}{|z|^n}\,.
$$
Changing to polar coordinates we get
$$
II_1 \le C\, \int_0^{2\,|x-y|} \sigma(\{\xi \in S^{n-1} : r \,\xi
\in S \})\, \frac{dr}{r}\,.
$$
By \eqref{eq7}
\begin{equation}\label{eq12}
\sigma(\{\xi \in S^{n-1} : r \,\xi \in S \}) \le C\,r^\ep \,
\end{equation}
and hence
$$
II_1 \le C\, |x-y|^\ep \,.
$$

One estimates  $II_2$ likewise, so we turn our attention to $III$\,.
The method is similar to what we have done before with other terms:
in the intersection of the domains of integration of the two
integrals in $III$ we apply a gradient estimate and in the
complement, which we split in four terms, we resort to the smallness
of the resulting domain of integration. Performing the plan just
sketched we get
\begin{equation*}
\begin{split}
III & \le \left|\int_{S \cap B(x,\R) \cap B^c(x,2\,|x-y|)\cap
B(y,\R) \cap B^c(y,2\,|x-y|)} \left(K(x-z) - K(y-z)\right) \,
dz\right|\\*[5pt] &\quad + \left|\int_{S \cap B(x,\R) \cap
B^c(x,2\,|x-y|)\cap B^c(y,\R)} K(x-z)\,dz \right|\\*[5pt] &\quad +
\left|\int_{S \cap B(x,\R) \cap B^c(x,2\,|x-y|) \cap B(y,2\,|x-y|)}
K(x-z)\,dz \right|\\*[5pt] &\quad + \left|\int_{S \cap B(y,\R) \cap
B^c(y,2\,|x-y|)\cap B^c(x,\R)} K(x-z)\,dz \right|\\*[5pt] &\quad +
\left|\int_{S \cap B(y,\R) \cap B^c(y,2\,|x-y|) \cap B(x,2\,|x-y|)}
K(x-z)\,dz \right|\\*[5pt] & =  III_1+ III_2+ III_3 +III_4+ III_5
\,.
\end{split}
\end{equation*}

By a gradient estimate $III_1$ is not greater than
$$
\left| \int_{S \cap B(x,\R) \cap B^c(x,2\,|x-y|)} C\,
\frac{|x-y|}{|x-z|^{n+1}}\,dz \right| \le C\,|x-y|\,
\int_{2\,|x-y|}^{\R} r^{-2+\ep }\,dr = C\, |x-y|^\ep \,,
$$
where \eqref{eq12} has been used in the first inequality\,.

The terms $III_2$ and $III_4$ are estimated in the same way. For
instance, for $III_2$  we have
\begin{equation*}
\begin{split}
III_2 & \le C\,\int_{B(x,\R) \cap B^c(y,\R)}
\frac{1}{|x-z|^n}\,dz\\*[7pt] & \le \frac{C}{\R^n}\, |B(x,\R) \cap
B^c(y,\R)|\\*[7pt] & \le  \frac{C}{\R}\, |x-y|\,.
\end{split}
\end{equation*}

The terms $III_3$ nd $III_5$ are also estimated in the same way. For
instance, for~\^{E}$III_3$ we have
$$
III_3 \le \int_{S \cap  B^c(x,2\,|x-y|) \cap B(x,3\,|x-y|)}
\frac{C}{|x-z|^n}\,dz \,.
$$
Since $x \in \Gamma$\,, for $z \in S \cap B(x,3\,|x-y|)$\,, we get
by Lemma 2
$$
|z| \le 6 \,(1+\sqrt{2})\, |x-y| \le 18 \, |x-y|\,,
$$
and so, making use of \eqref{eq12}\,,
\begin{equation*}
\begin{split}
III_3 &\le C\, \int_{S \cap B(0, 18\, |x-y|)} \frac{dz}{|z|^n} \le
C\,|x-y|^\ep \,.
\end{split}
\end{equation*}

\vspace*{7pt}

\noindent {\it Case 2:} $x$ and $y$ are in tangential position, that
is, $y \in \Omega \setminus \Gamma\,.$ We intend to perform a
reduction to the non-tangential case. With this in mind take the
point $p$ in $\fr$ nearest to $y$ and let $\overrightarrow{N}$ be
the inner unit normal vector to $\fr$ at the point $p$\,. Consider
the ray $y+t \overrightarrow{N }\,,\,\, t > 0$\,. See the Figure 4.

\begin{figure}[ht]
\begin{center}
\includegraphics{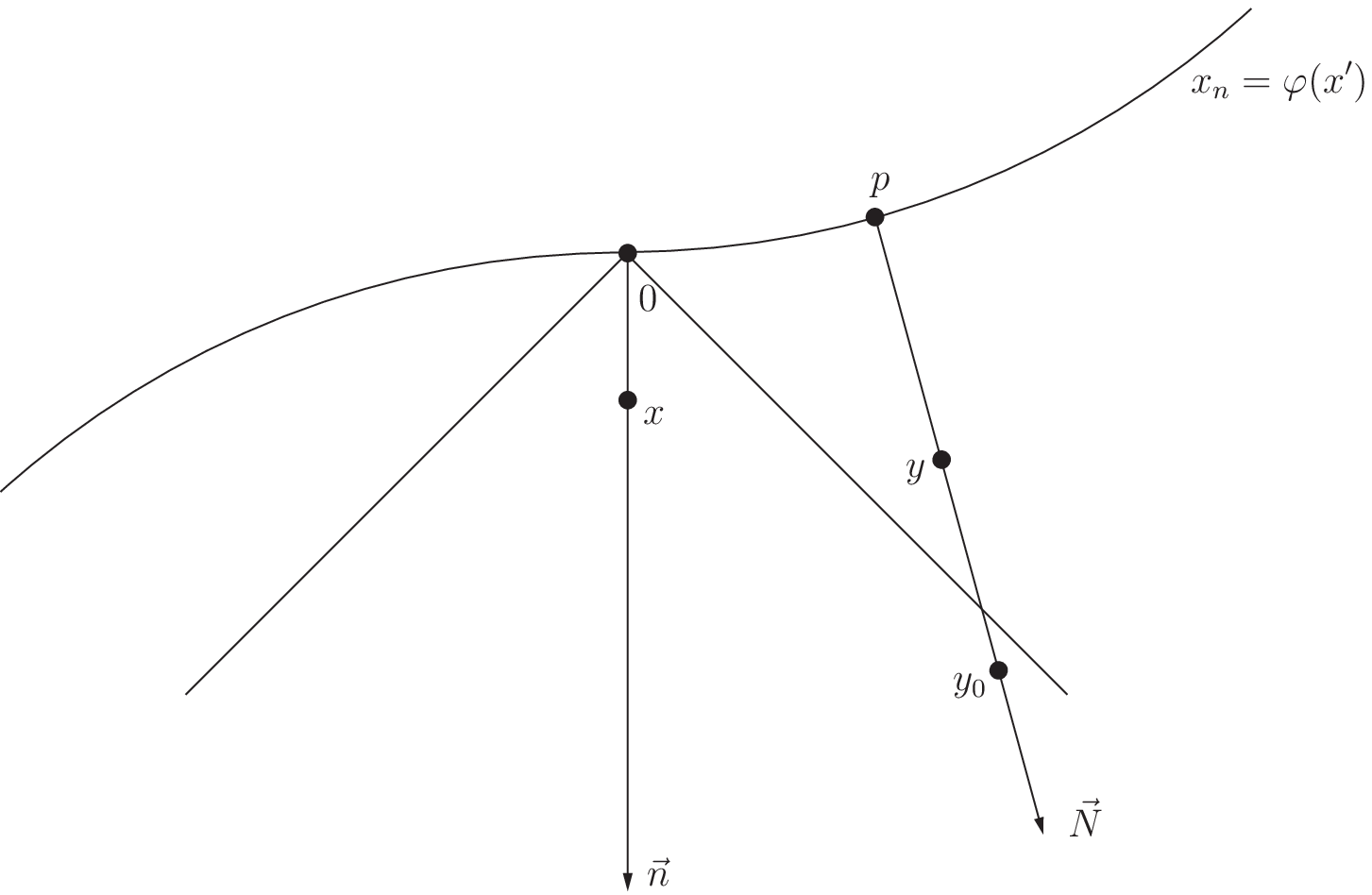}\\*[7pt]
Figure 4
\end{center}
\end{figure}

\noindent
 The condition on $t$ for $y+t \overrightarrow{N} \in
\Gamma$\, is
\begin{equation}\label{eq13}
(y+t \overrightarrow{N},\overrightarrow{n}) \geq
\frac{1}{\sqrt{2}}\, |y+t \overrightarrow{N}|\,.
\end{equation}
We clearly have
$$
(y+t \overrightarrow{N},\overrightarrow{n}) \geq
t\,(\overrightarrow{N}, \overrightarrow{n})- |y|
$$
and
$$
|y+t \overrightarrow{N}| \leq |y|+t \,.
$$
A sufficient condition for \eqref{eq13} is then
$$
t \geq \frac{(1+\frac{1}{\sqrt{2}})\,|y|}{(\overrightarrow{N},
\overrightarrow{n})- \frac{1}{\sqrt{2}}}\,.
$$
If $\R$ is small enough, then $(\overrightarrow{N},
\overrightarrow{n}) \geq \displaystyle\frac{3}{4}
\frac{1}{\sqrt{2}}$\,, and thus we obtain a simpler sufficient
condition for \eqref{eq13}\, namely,
$$
t \geq 3 \,(1+\sqrt{2})\,|y| \equiv t_0\,.
$$
Set $y_0 = y+t_0 \,\overrightarrow{N}$\,, so that $y_0 \in
\Gamma$\,. The reduction will be completed if we show that
\begin{equation}\label{eq14}
|y-y_0| \le C \, |x-y|\,,
\end{equation}
because $x$ and $y_0$ on one hand, and $y$ and $y_0$ on the other,
are in non-tangential position.

Clearly $|y-y_0|= t_0 = C\,|y|$\,. Since $y \in \Omega \setminus
\Gamma$\,, $|(y,\overrightarrow{n})| < \frac{1}{\sqrt{2}}\,|y|$\,,
and hence
$$
|y| \le |(y,\overrightarrow{n})|+ |y'| \le
\frac{1}{\sqrt{2}}\,|y|+|y'|\,,
$$
which yields $ (1-\frac{1}{\sqrt{2}})\,|y| \le |y'|\,.$ Therefore
\begin{equation*}
\begin{split}
|x-y| & = \left|
|x|\,\overrightarrow{n}-(y,\overrightarrow{n})\overrightarrow{n}-y'\right|
\\*[5pt] & \geq  |y'| \geq
\left(1-\frac{1}{\sqrt{2}}\right)\,|y|\\*[5pt] & = C\, |y-y_0|\,,
\end{split}
\end{equation*}
which is \eqref{eq14}\,.

This completes the proof that $T(\chi_\Omega) \in \Lipo$.

We now proceed to prove that for an arbitrary $f \in \Lipo$ one has
that $ T(f) \in \Lipo$. Recall that we already know that $T(f) \in
L^\infty(\Omega)$\,(see \eqref{eq8})\,. To estimate the semi-norm
$\sigma_\ep(T(f))$ we start with the obvious decomposition of
$T(f)$, namely,
$$
T(f)(x)= \int_\Omega (f(y)-f(x))\, K(x-y)\,dy + f(x)
T(\chi_\Omega)(x) \,,
$$
so that only the first term, which we denote by $S(f)(x)$\,, is
still a problem. Let $A$ stand for $\Omega \cap B(x_1, 2 |x_1-x_2|)$
and set $B= \Omega \setminus A$\,. Then
\begin{equation*}
\begin{split}
S(f)(x_1)\!-\!S(f)(x_2) &  = \int_A \left[ (f(y)-f(x_1)) \, K(x_1-y)
- (f(y)-f(x_2) ) \, K(x_2-y)\right]\,dy \\*[7pt] &\quad +\! \int_{B}
\left[ (f(y)-f(x_1)) \, K(x_1-y) \!-\! (f(y)\!-\!f(x_2)) \,
K(x_2-y)\right]\,dy\\*[7pt] & = I+II\,.
\end{split}
\end{equation*}
Set $A' =\Omega \cap B(x_2, 3 |x_1-x_2|)$. Clearly $A \subset A'$\,.
The term $I$ is easy to estimate as indicated below.
\begin{equation*}
\begin{split}
|I| & \leq C\, \left(\int_A \frac{|f(y)-f(x_1)|}{|y-x_1|^{n}}\,dy +
\int_{A'} \frac{|f(y)-f(x_2)|}{|y-x_2|^{n}}\,dy \right) \\*[7pt] &
\le C\,\|f\|_\ep \, \int_0^{3\,|x_1-x_2|} r^{-1+\ep}\,dr \\*[7pt] &
= C\,\|f\|_\ep \,|x_1-x_2|^\ep \,.
\end{split}
\end{equation*}
For the term $II$ we have
\begin{equation*}
\begin{split}
II & = \int_B (f(y)-f(x_2))\,\left(K(x_1-y)-K(x_2-y)\right)\,dy
\\*[7pt] &\quad + (f(x_2)-f(x_1))\, \int_{B} K(x_1-y)\,dy \\*[7pt] &
= III+IV\,.
\end{split}
\end{equation*}
On one hand, by a gradient estimate we get
\begin{equation*}
\begin{split}
|III| & \le \|f\|_\ep \,|x_1-x_2|\, \int_B
\frac{1}{|y-x_1|^{n+1-\ep}} \,dy \\*[7pt] & = \|f\|_\ep
\,|x_1-x_2|\, \int_{2\,|x_1-x_2|}^\infty r^{-2+\ep}\,dr \\*[7pt] & =
C\,\|f\|_\ep \,|x_1-x_2|^\ep\,,
\end{split}
\end{equation*}
and on the other hand, we clearly have by \eqref{eq8}
$$
|IV| \le \|f\|_\ep \,|x_1-x_2|^\ep \, T^*(\chi_\Omega)(x_1) \le
C\,\|f\|_\ep \,|x_1-x_2|^\ep\,.
$$

Finally, one can check without pain that the arguments above may be
adapted to yield the boundedness of $T$ as a map from $\Lipo$ into
$\operatorname{Lip}(\ep,\Omega^c)$.
\end{proof}

\section{Estimates for Commutators}

In this section we consider the commutator between the smooth
homogeneous even Calder\'{o}n-Zygmund operator $T$ (see \eqref{eq6bis})
and the multiplication operator by a function $a \in
\operatorname{Lip}(\alpha,\Omega)\,, 0<\alpha<1 \,,$
\begin{equation}\label{eq15}
[T,a](f)(x)= \int_\Omega \left(
a(x)-a(y)\right)\,K(x-y)\,f(y)\,dy\,,\quad x \in \Omega\,,
\end{equation} where $K(x)$ is the kernel of $T$ and $f \in
\operatorname{Lip}(\beta,\Omega)\,, \,0<\beta<1 $\,. As in the
previous section, $\Omega$ is a bounded domain in $\Rn$ with smooth
boundary of class $C^{1+\ep}$\,, $0<\ep <1\,.$

\begin{lemma}\label{lemma 3}
For $0<\alpha<1$ and $0 < \beta \leq \ep$ we have the estimate
\begin{equation}\label{eq16}
\| [T,a] (f)\|_\alpha \le C\, \sigma_\alpha(a)\, \|f\|_\beta\,,\quad
f \in \operatorname{Lip}(\beta,\Omega)\,,
\end{equation}
where $C$ is a constant depending only on $n$, $\Omega$, $\ep$,
$\alpha$ and $\beta$\,.
\end{lemma}
Recall that for $0<\alpha<1$\,,
$$
\|g \|_\alpha = \|g\|_\infty + \sigma_\alpha(g)\,,
$$
where $\|g\|_\infty$ is the supremum norm of $g$ on $\Omega$ and
$$
\sigma_\alpha(g)= \sup \left\{ \frac{|g(x)-g(y)|}{|x-y|^{\alpha}} :
x,y \in \Omega , \,x \neq y \right\}\,.
$$

A consequence of the preceding lemma is that if $\beta \le \ep$ and
$\beta < \alpha$ then the commutator $[T,a]$ is compact as an
operator from $\operatorname{Lip}(\beta,\Omega)$ into itself\,. This
follows from the fact that each ball of
$\operatorname{Lip}(\alpha,\Omega)$ is relatively compact in
$\operatorname{Lip}(\beta,\Omega)$\,(\cite[Corollary 3.3,
p.~154]{J}). The Lemma is applied to the Beurling transform and the
function $a=\mu$\,. Then $\alpha = \ep$ and $\beta =\ep'$ \,, where
$\ep'$ is any number with $0<\ep'<\ep$\,.

\begin{proof}
We first estimate $\| [T,a] (f)\|_\infty$\,. For each $x \in
\Omega$,
\begin{equation*}
\begin{split}
| [T,a] (f)(x)| &  \le C\,\sigma_\alpha(a)\, \|f\|_\infty \,
\int_{\Omega} |x-y|^{-n+\alpha} \,dy\\*[7pt] & \le
C\,\sigma_\alpha(a)\, \|f\|_\infty \, \int_0^{d} r^{-1+\alpha}\,dr
\\*[7pt] & = C\,\sigma_\alpha(a)\, \|f\|_\infty \, d^\alpha \,,
\end{split}
\end{equation*}
where $d$ is the diameter of $\Omega$\,.

We turn now to the more difficult task of estimating
$\sigma_\alpha([T,a](f))$\,.
 Fix
$x_1$ and~$x_2$ in $\Omega$\,. Then
\begin{equation*}
\begin{split}
|[T,a](f)(x_1)\!-\![T,a] (f)(x_2)| & \le
|a(x_1)-a(x_2)|\left|\int_\Omega K(x_1-y)f(y)\,dy)\right| \\*[7pt]
&\quad +\!\left|\int_\Omega
(a(x_2)-a(y))\,(K(x_1-y)-\!\!K(x_2-y))\,f(y)\,dy \right|\\*[7pt] & =
I + II\,,
\end{split}
\end{equation*}
and clearly, by the Main Lemma,
$$
I \le C\,|x_1-x_2|^\alpha\,\sigma_\alpha(a)\, \|T(f)\|_\infty \le
C\,|x_1-x_2|^\alpha\,\sigma_\alpha(a)\, \|f\|_\beta \,.
$$
To estimate $II$ we introduce the sets
$$
A= \{y \in \Omega : |y-x_1|> 2\,|x_1-x_2|\}
$$
and
$$
B= \{y \in \Omega : |y-x_1| \leq 2\,|x_1-x_2|\}\,.
$$
Notice that $|y-x_2|> |x_1-x_2|\,,\,\, y \in A$ and $|y-x_2|\leq 3\,
|x_1-x_2|\,,\,\, y \in B$\,. Let~$II_A$ (respectively $II_B$) denote
the absolute value of the integral in $II$ with domain of
integration restricted to $ A$ (respectively to $B $).

By a gradient estimate
\begin{equation*}
\begin{split}
 II_A & \le \int_A |a(x_2)-a(y)|\,\frac{|x_1-x_2|}{|x_2 -
y|^{n+1}}\,|f(y)|\,dy \\*[7pt] & \le
C\,|x_1-x_2|\,\sigma_\alpha(a)\, \|f\|_\infty \int_{A}
|x_2-y|^{-(n+1)+\alpha}\,dy\\*[7pt] & \le
C\,|x_1-x_2|\,\sigma_\alpha(a)\, \|f\|_\infty
\int_{|x_1-x_2|}^\infty r^{-2+\alpha}\,dr\\*[7pt] & \le
C\,\sigma_\alpha(a)\, \|f\|_\infty \,|x_1-x_2|^\alpha\,.
\end{split}
\end{equation*}
For the term $II_B$ we have
\begin{equation*}
\begin{split}
 II_B & \le \left|\int_B  (a(x_2)-a(y))\,(K(x_1-y)\,f(y)\,dy \right|\\*[7pt]
&\quad + \left|\int_B  (a(x_2)-a(y))\,(K(x_2-y)\,f(y)\,dy \right|
\\*[7pt] & = III+IV\,,
\end{split}
\end{equation*}
and $IV$ can be estimated directly as follows
\begin{equation*}
\begin{split}
 IV & \le   \sigma_\alpha(a)\, \|f\|_\infty \, \int_B \frac{|x_2-y|^\alpha}{|x_2 -
y|^{n}}\,dy \\*[7pt] & \le   \sigma_\alpha(a)\, \|f\|_\infty \,
\int_0^{3\,|x_1-x_2|} r^{-1+\alpha}\,dr \\*[7pt] & = C\,
\sigma_\alpha(a)\, \|f\|_\infty \,|x_1-x_2|^\alpha\,.
\end{split}
\end{equation*}
The term $III$ needs an additional manoeuvre, which consists is
bringing back $a(x_1)$:
\begin{equation*}
\begin{split}
 III & \le \left|\int_B  (a(x_1)-a(y))\,(K(x_1-y)\,f(y)\,dy \right| \\*[7pt]
&\quad +  |a(x_2)-a(x_1)|\left|\int_B K(x_1-y)f(y)\,dy)\right|
\\*[7pt] & = IV' + V\,,
\end{split}
\end{equation*}
and $IV'$ can be treated as $IV$\,. Now
$$
\int_B K(x_1-y)f(y)\,dy = \int_\Omega K(x_1-y)f(y)\,dy -
\int_{\Omega \cap B^c(x_1,\,2\,|x_1-x_2|)} K(x_1-y)f(y)\,dy
$$
and thus, by \eqref{eq8} ,
$$
\left|\int_B K(x_1-y)f(y)\,dy\right| \le 2\,T^*(f)(x_1) \le C\,
\|f\|_\beta \,.
$$
Therefore
$$
V \le C\,   \sigma_\alpha(a)\, \|f\|_\beta \,|x_1-x_2|^\alpha\,.
$$
\end{proof}

\section{Relationship between \boldmath$B_\Omega^n$ and \boldmath$B^n$}

Recall that if $B$ is the Beurling transform then $B_\Omega(f):=
B(f)\,\chi_\Omega\,.$ The main goal of this section is to prove
the following result.

\begin{teorema}\label{th1}
Let $\Omega \subset \C$ be a  bounded domain with boundary of class
$C^{1+\ep}$, $0<\ep <1$. Then, for each positive integer $n$\,, we
have
$$
B_\Omega^n(f)(z) = B^n(f)(z) \chi_\Omega(z) + K_n(f)(z)\,,
$$
where $K_n$ is a compact operator from $\Lip$ into itself,
$0<\ep'<\ep$\,.
\end{teorema}

\begin{proof} For $n\geq 2$, we obtain, proceeding by induction,
\begin{equation*}
\begin{split}
B_\Omega^n(f)& = B(B_\Omega^{n-1}(f))\,\chi_\Omega \\*[5pt] &=
B\left(B^{n-1}(f)\,\chi_\Omega+K_{n-1}(f)\right)\,\chi_\Omega
\\*[5pt] &= B\left(B^{n-1}(f)-
B^{n-1}(f)\,\chi_{\Omega^c}+K_{n-1}(f)\right)\,\chi_\Omega\\*[5pt] &
= B^n(f)\,\chi_\Omega - B(B^{n-1}(f)\,\chi_{\Omega^c})\,\chi_\Omega
+B(K_{n-1}(f))\,\chi_\Omega\,.
\end{split}
\end{equation*}
It is then enough to prove that, for $n\geq1$, the operator
$$
B(B^{n}(f)\,\chi_{\Omega^c})\,\chi_{\Omega}
$$
is compact from $\Lip$ into itself.

Let $dA$ stand for area measure in the plane and take a function $f
\in L^\infty(\Omega)$. Then, for $z \in \Omega$\,,
\begin{equation*}
\begin{split}
B(B^{n}(f)\,\chi_{\Omega^c})(z) & = -\frac{1}{\pi} \int_{\Omega^c}
\frac{B^{n}(f)(w)}{(z-w)^2}\,dA(w) \\*[7pt] & =- \frac{1}{\pi}
\int_{\Omega^c} \frac{1}{(z-w)^2} \frac{(-1)^n n}{\pi}\,\int_\Omega
\frac{(\overline{w-\zeta})^{n-1}}{(w-\zeta)^{n+1}}\,f(\zeta)\,dA(\zeta)
\,dA(w) \\*[7pt] & = c_n\, \int_\Omega
K(z,\zeta)\,f(\zeta)\,dA(\zeta)\,,
\end{split}
\end{equation*}
where
$$
K(z,\zeta) = K_n(z,\zeta):= \int_{\Omega^c} \frac{1}{(z-w)^2}
\frac{n \,(\overline{w-\zeta})^{n-1}}{(w-\zeta)^{n+1}}\,dA(w)
$$
and $c_n = \frac{(-1)^{n+1} }{\pi^2}\,.$

Notice that if $\Omega$ is a disc, say the unit disc, then
$K(z,\zeta)=0,\,\,z,\zeta \in \Omega\,.$ To see this readily, apply
Green-Stokes' Theorem to the complement of the unit disc to obtain
$$
K(z,\zeta)= \frac{-1}{2\,\imath} \int_{\fr} \frac{1}{(z-w)^2}
\frac{(\overline{w-\zeta})^{n}}{(w-\zeta)^{n+1}}\,dw \,.
$$
Expand $(\overline{w-\zeta})^{n}$ by Newton's formula and then use
$\overline{w} = \frac{1}{w}\,, \,\, |w|=1 \,.$  Thus $K(z,\zeta)$~is
a finite sum of integrals over the unit cercle of rational functions
with all poles in the open unit disc. Hence each of these integrals
is zero.

We claim that if $\Omega$ is not a disc, then  the operator
$$
P(f)(z)= \int_\Omega K(z,\zeta)\,f(\zeta)\,dA(\zeta),\quad
z\in\Omega\,,
$$
which may be non-zero,  is a smoothing operator. By this we mean
that
\begin{equation}\label{eq17}
 \|P(f)\|_\alpha \le C\, \|f\|_\infty\,,\quad 0<
\alpha<\ep\,,
\end{equation}
where $C$ depends only on $\alpha$\,, $\ep$ and $\Omega$\,.

Of course \eqref{eq17} completes the proof of  Theorem 1, because
then $P$ maps the unit ball of $\Lip$ into a ball of
$\operatorname{Lip}(\alpha,\Omega)$, for $ \alpha < \ep$\,, which is
relatively compact in $\Lip$ provided $\ep'< \alpha$.

Our next goal is to show that \eqref{eq17} is a consequence of the
properties of the kernel~$K(z,\zeta)$ described in the following
lemma.

\begin{lemma}\label{lemma4}
The kernel $K(z,\zeta)$ satisfies the following
\begin{alignat*}{1}
&\text{(i)}\hskip3cm |K(z,\zeta)| \le C\,
\frac{1}{|z-\zeta|^{2-\ep}}\,,\quad z\,,\,\zeta \in
\Omega\,.\\*[7pt] &\text{(ii)}\hskip1cm|K(z_1,\zeta)-K(z_2,\zeta)|
\le C \, \frac{|z_1-z_2|^\varepsilon}{|\zeta-z_1|^2}\,,\quad
z_1,\,z_2 \in \Omega\,,\quad |\zeta-z_1| \geq 2\,|z_1-z_2|\,.
\end{alignat*}
\end{lemma}

Before discussing the proof of Lemma 6 we show how it yields
\eqref{eq17}\,.

We first prove that $P(f)$ is bounded on $\Omega$. Denoting by $d$
the diameter of $\Omega$, we obtain, by Lemma 6 (i)\,,
\begin{equation*}
\begin{split}
|P(f)(z)| &\le \int_ \Omega |K(z,\zeta)|
|f(\zeta)|\,dA(\zeta)\\*[7pt] & \le C\, \|f\|_\infty\, \int_ \Omega
\frac{dA(\zeta)}{|z-\zeta|^{2-\ep}}\\*[7pt] & \le C\, \|f\|_\infty\,
\int_0^d r^{-1+\ep}\,dr \\*[7pt] & = C\,d^\ep \, \|f\|_\infty \, .
\end{split}
\end{equation*}

Next we claim that
\begin{equation}\label{eq18}
|P(f)(z_1)-P(f)(z_2)| \le C\,|z_{1}-z_{2}|^\ep \,(1+\log
\frac{d}{|z_{1}-z_{2}|})\,\|f\|_\infty\,, \quad z_1\,,\,z_2 \in
\Omega \,.
\end{equation}
Clearly \eqref{eq17} follows from \eqref{eq18}\,. To prove
\eqref{eq18} take $z_1\,,\,z_2 \in \Omega$\,.  Define $A= \{\zeta
\in \Omega : |z_1-\zeta| < 2\,|z_1-z_2|\}$ and $B=\Omega \setminus
A$\,. Therefore
\begin{equation*}
\begin{split}
|P(f)(z_1)-P(f)(z_2)| & \le \int_A |K(z_1,\zeta)|
|f(\zeta)|\,dA(\zeta)\\*[7pt] &\quad +\int_A |K(z_2,\zeta)|
|f(\zeta)|\,dA(\zeta)\\*[7pt] &\quad+ \int_B
|K(z_1,\zeta)-K(z_2,\zeta)| |f(\zeta)|\,dA(\zeta)\\*[7pt] & =
I+II+III \,.
\end{split}
\end{equation*}
Applying Lemma 6 (i), the terms $I$ and $II$ can be estimated by
$$
C\,\|f\|_\infty\, \int_0^{3\,|z_1-z_2|} r^{-1+\ep} \,dr \le
C\,|z_1-z_2|^\ep \,\|f\|_\infty\,.
$$
Applying Lemma 6 (ii)\,, the term $III$ can be estimated by
\begin{equation*}
\begin{split}
III & \le C\,\|f\|_\infty\, |z_1-z_2|^\varepsilon \left(\int_B
\frac{dA(\zeta)}{|\zeta -z_1|^{2}}\right) \\*[7pt] & \le
C\,\|f\|_\infty\, |z_1-z_2|^\varepsilon \, \int_{2\,|z_1-z_2|}^d
\frac{dr}{r} \\*[7pt] & = C\, \|f\|_\infty\, |z_1-z_2|^\varepsilon
\,\log{\frac{d}{2\,|z_1-z_2|}}\,,
\end{split}
\end{equation*}
which completes the proof of Theorem 1.
\end{proof}

\begin{proof}[Proof of Lemma~\ref{lemma4}]

For each $\zeta\in\Omega$ consider the Cauchy integral of the
function $(\overline{w-\zeta})^n$ on~$\partial\Omega$, that is,
$$
H_{\zeta}(w)=\frac{1}{2\pi
i}\int_{\partial\Omega}\frac{(\overline{t-\zeta})^n}{t-w}\,dt,\quad
w\in\mathbb{C}\setminus \partial\Omega.
$$
For $w\in\partial\Omega$ let $H_{\zeta}(w)$ be the non-tangential
limit of $H_{\zeta}$ from~$\Omega$, that is, the limit of
$H_{\zeta}(w')$ as $w'\in\Omega$ tends to~$w$ non-tangentially.
Similarly, denote by $H^c_{\zeta}(w)$ the non-tangential limit of
$H_{\zeta}$ from $\mathbb{C}\setminus\Omega$. These limits exist
a.e.\ on $\partial\Omega$ with respect to arc-length and one has the
Plemelj formula (e.g. \cite[p.~143]{Ve})
$$
(\overline{w-\zeta})^n=H_{\zeta}(w)-H^c_{\zeta}(w),\quad w\text{
a.e.\ on }\partial\Omega.
$$
Indeed, it can be shown that $H_{\zeta}$ is of
class~$C^{1+\varepsilon}$ in~$\Omega$ and in $\mathbb{C}\setminus
\overline{\Omega}$, so that the above limits exist everywhere
on~$\partial\Omega$ and without the non-tangential approach
restriction. We do not need, however, such fact.

Applying the Green-Stokes Theorem to the form
$$
\frac{(\overline{w-\zeta})^n+H^c_{\zeta}(w)}{(z-w)^2(w-\zeta)^{n+1}}\,dw
$$
and the domain $\Omega^c$, we get
$$
K(z,\zeta)=-\int_{\partial\Omega}\frac{H_{\zeta}(w)}{(z-w)^2(w-\zeta)^{n+1}}\,dw,
$$
which by the Residue Theorem is
$$
-2\pi i\left\{
\left.\frac{d}{dw}\frac{H_{\zeta}(w)}{(w-\zeta)^{n+1}}\right|_{w=z}+\frac{1}{n!}
\left.\frac{d^n}{dw^n}\frac{H_{\zeta}(w)}{(w-z)^2}\right|_{w=\zeta}\right\}.
$$
A straightforward computation of the residues yields
\begin{multline}\label{eq22}
K(z,\zeta)=-2\pi i\left\{
\frac{H'_{\zeta}(z)}{(z-\zeta)^{n+1}}-(n+1)
\frac{H_{\zeta}(z)}{(z-\zeta)^{n+2}}\right.\\*[5pt]
\left.+\sum^n_{\ell=0}(-1)^{n-\ell}\frac{(n-\ell+1)}{\ell!}
\frac{d^\ell}{d\zeta^\ell}H_{\zeta}(\zeta)\frac{1}{(\zeta-z)^{n+2-\ell}}.
\right\}.
\end{multline}
In the expression above for the kernel $K(z,\zeta)$ one may divine
the presence of non-obvious cancellation properties (consider the
case $n=1$). The strategy to unravel them is to bring into the
scene the function
 $$h(z)=2\pi i\,
H_{z}(z)=\displaystyle\int_{\partial\Omega}\frac{(\overline{t-z})^n}{t-z}\,dt\,,$$
and express $K(z,\zeta)$ in terms of $h$ and its derivatives.
Taylor's expansions of $h$ and its derivatives will then help in
understanding cancellations. The derivatives of $h$ are given by
\begin{equation}\label{eq23}
\frac{\partial^{\ell}}{\partial
z^{\ell}}\frac{\partial^k}{\partial\overline{z}^k}h(z)=(-1)^k\frac{\ell!\,n!}{(n-k)!}
\int_{\partial\Omega}\frac{(\overline{t-z})^{n-k}}{(t-z)^{1+\ell}}\,dt.
\end{equation}
On the other hand, by the binomial formula,
\begin{equation*}
\begin{split}
2\pi i\, H_{\zeta}(z)&=
\int_{\partial\Omega}\frac{(\overline{t-\zeta})^n}{t-z}\,dt\\*[5pt]
&=\sum^n_{\ell=0}\binom{n}{\ell} (\overline{z-\zeta})^{\ell}
\int_{\partial\Omega}
\frac{(\overline{t-z})^{n-\ell}}{t-z}\,dt\\*[5pt]
&=\sum^n_{\ell=0}\frac{(-1)^\ell}{\ell!}
\frac{\partial^{\ell}h}{\partial \overline{z}^\ell}(z)
(\overline{z-\zeta})^{\ell}.
\end{split}
\end{equation*}
Differentiating the preceding identity with respect to $z$
\begin{equation}\label{eq24}
2\pi i\,
H'_{\zeta}(z)=\sum^n_{\ell=0}\frac{(-1)^{\ell}}{\ell!}\frac{\partial^{\ell+1}}{\partial\overline{z}^{\ell}\partial
z}h(z)(\overline{z-\zeta})^{\ell}.
\end{equation}
Therefore
\begin{equation}\label{eq25}
\begin{split}
-(\zeta-z)^{n+2}
 K(z,\zeta)&=
\sum^n_{\ell=0} \frac{(-1)^{\ell}}{\ell!}
\frac{\partial^{\ell+1}}{\partial\overline{z}^{\ell}\partial
z}h(z)(\overline{z-\zeta})^{\ell}(z-\zeta)\\*[5pt]
&\quad-(n+1)\sum^n_{\ell=0}\frac{(-1)^{\ell}}{\ell!}
\frac{\partial^{\ell}h}{\partial \overline{z}^{\ell}}(z)
(\overline{z-\zeta})^{\ell}\\*[5pt] &\quad+
\sum^n_{\ell=0}\frac{(n+1-\ell)}{\ell!}
\frac{\partial^{\ell}h(\zeta)}{\partial \zeta^\ell}(z-\zeta)^{\ell}.
\end{split}
\end{equation}
In each of the terms of the last sum it will be convenient to write
a Taylor expansion of
$\displaystyle\frac{\partial^{\ell}h(\zeta)}{\partial \zeta^{\ell}}$
up to order~$n-\ell$ around the point~$z$. Doing so we obtain
\begin{equation*}
\begin{split}
(z-\zeta)^{n+2}K(z,\zeta)&= \sum^n_{\ell=0} \frac{1}{\ell!}
\frac{\partial^{\ell}}{\partial
\overline{z}^{\ell}}\frac{\partial}{\partial z}
h(z)(\overline{\zeta-z})^{\ell} (\zeta-z)\\*[5pt] &\quad
+(n+1)\sum^n_{\ell=0}\frac{1}{\ell!}
\frac{\partial^{\ell}h}{\partial \overline{z}^{\ell}} (z)
(\overline{\zeta-z})^{\ell}\\*[5pt] &\quad-\sum^n_{\ell=0}
\frac{n+1-\ell}{\ell!}\\*[5pt] &\qquad\times
\sum^{n-\ell}_{j=0}\sum^{j}_{k=0} \frac{(-1)^{\ell}}{k!\,(j-k)!}
\frac{\partial^{\ell+k}}{\partial z^{\ell+k}}
\frac{\partial^{j-k}}{\partial
\overline{z}^{j-k}}h(z)(\zeta-z)^{k+\ell}(\overline{\zeta-z})^{j-k}\\*[5pt]
&\quad+ R(z,\zeta)\equiv S(z,\zeta)+R(z,\zeta).
\end{split}
\end{equation*}

A cumbersome but easy computation shows now that
$$
S(z,\zeta)=0,\quad z,\zeta\in\Omega.
$$
 The most direct way to ascertain this is to check that the
coefficient of $S(z,\zeta)$ in the monomial
$(\zeta-z)^{m_0}\,(\overline{\zeta-z})^{p_0}$ vanishes for all
non-negative exponents $m_0$ and $p_0.$  For this we distinguish
four cases.

\vspace*{7pt} \noindent {\it Case 1:} Assume that $m_0 \ge 2.$
Only in the third sum may appear terms of this type and they must
cancel out by themselves.  This can be shown using the identities
$$
\sum^{m_{0}}_{\ell=0}\binom{m_{0}}{\ell}(-1)^{\ell}
=\sum^{m_{0}}_{\ell=0}\ell\binom{m_{0}}{\ell} (-1)^{\ell}=0.
$$

\vspace*{7pt} \noindent {\it Case 2:} Take $m_0=1$ and $0 \le p_0
\le n-1 .$  Two terms appear in the third sum and one in the
third, and they cancel.

\vspace*{7pt} \noindent {\it Case 3:} Take $m_0=1$ and $p_0 = n .$
There is only one term of this type, which corresponds to letting
$l=n$ in the first sum. To show that this term vanishes we resort
to \eqref{eq23} for $l=1$ and $k=n$ and then we apply Cauchy's
Theorem.

\vspace*{7pt} \noindent {\it Case 4:} Take $m_0=0$ and $0 \le p_0
\le n .$  One term in the second sum cancels with a term in the
third sum.

We turn now to the analysis of the kernel $K(z,\zeta).$ Since
$S(z,\zeta)$ vanishes identically we get
\begin{equation}\label{eq25bis}
-(\zeta-z)^{n+2}
 K(z,\zeta)= \sum^n_{\ell=0}\frac{(n+1-\ell)}{\ell!}
\, R_{n-l}(z,\zeta)\,(z-\zeta)^{\ell}\,,
\end{equation}
where $R_{n-l}(z,\zeta)$ is the remainder of the Taylor expansion of
$\displaystyle\frac{\partial^{\ell}h(\zeta)}{\partial \zeta^{\ell}}$
up to order~$n-\ell$ around the point~$z$.

 A key fact in the present proof is that the remainder $R_{n-l}(z,\zeta)$
is $O(|z-\zeta|^{n-l+\varepsilon})$, because the $n$-th order
derivatives of $h(z)$ are in $\operatorname{Lip}
(\varepsilon,\Omega)$. To show this we resort to \eqref{eq23} to get
$$
\frac{\partial^k}{\partial
z^k}\frac{\partial^{n-k}}{\partial\overline{z}^{n-k}}h(z)=(-1)^{n-k}\,n!
\int_{\partial\Omega}\frac{(\overline{t-z})^k}{(t-z)^{1+k}}\,dt.
$$
If $k=0$, then the above expression is
$$
(-1)^n n!\int_{\partial\Omega}\frac{dt}{(t-z)}=(-1)^n n!\,2\pi i.
$$
If $1\le k\le n$, then we obtain, by Green-Stokes and for some
constant~$c_{n,k}$,
$$
(-1)^{n-k}n!\,
k\,2i\int_{\Omega}\frac{(\overline{t-z})^{k-1}}{(t-z)^{k+1}}\,dA(t)=
c_{n,k}\, B^k(\chi_{\Omega})(z),
$$
which is in $\operatorname{Lip}(\varepsilon,\Omega)$ owing to the
Main Lemma. Here $B^k$ is the $k$-th iteration of the Beurling
transform.

Part (i) of the lemma is a straightforward consequence of
\eqref{eq25bis} and the size estimate on the remainder
$R_{n-l}(z,\zeta)$ we have just proved.

We are left with part (ii). Take points $z_1,z_2$ and $\zeta$ in
$\Omega$ with $|\zeta-z_1| \geq 2\,|z_1-z_2|\,.$  From
\eqref{eq25bis} we obtain
\begin{equation*}
K(z_1,\zeta)-K(z_2,\zeta) = \frac{(-1)^{n+1}}{2 \pi
\imath}\,\sum^n_{\ell=0}\frac{(n+1-\ell)}{\ell!} \, \left(
\frac{R_{n-l}(z_1,\zeta)} {(z_1-\zeta)^{n+2-l}} -
\frac{R_{n-l}(z_2,\zeta)} {(z_2-\zeta)^{n+2-l}}\right)
\end{equation*}
Add and subtract $R_{n-l}(z_1,\zeta)$ in the numerator of the second
fraction above to get

\begin{equation*}
K(z_1,\zeta)-K(z_2,\zeta) = I + II\,,
\end{equation*}
where
\begin{equation*}
I = \frac{(-1)^{n+1}}{2 \pi
\imath}\,\sum^n_{\ell=0}\frac{(n+1-\ell)}{\ell!} \,
R_{n-l}(z_1,\zeta) \, \left(\frac{1}{(z_1-\zeta)^{n+2-l}}-
\frac{1}{(z_2-\zeta)^{n+2-l}} \right)
\end{equation*}
and
\begin{equation*}
II = \frac{(-1)^{n+1}}{2 \pi
\imath}\,\sum^n_{\ell=0}\frac{(n+1-\ell)}{\ell!}\,\frac{R_{n-l}(z_1,\zeta)-R_{n-l}(z_2,\zeta)}{(z_2-\zeta)^{n+2-l}}\,.
\end{equation*}
Controlling $I$ is easy via an obvious gradient estimate, which
yields
\begin{equation*}
\begin{split}
I & \le C\, |z_1-\zeta|^{n-l+\varepsilon}\,
\frac{|z_1-z_2|}{|z_1-\zeta|^{n+3-l}}
\\*[7pt] & =
C\,\frac{|z_1-z_2|}{|z_1-\zeta|^{3-\varepsilon}}
\\*[7pt] &  \le C\, \frac{|z_1-z_2|^\varepsilon}{|z_1-\zeta|^2}\,.
\end{split}
\end{equation*}

To estimate the term $II$ we need a sublemma.
\begin{sublemma}\label{sublema} We have the identity
\begin{equation*}
\begin{split}
 R_{n-l}(z_1,\zeta)-R_{n-l}(z_2,\zeta) & = \sum_{j+k = n-l}
c_{j,k}\,\left(
B^{l+j}(\chi_\Omega)(z_1)-B^{l+j}(\chi_\Omega)(z_2)\right)\,(\zeta-z_2)^j\,(\overline{\zeta-z_2})^k
\\ & + O(|z_1 - z_2|^{1+\varepsilon} \, |\zeta-z_2|^{n-l-1})\,.
\end{split}
\end{equation*}
\end{sublemma}
\vspace{0.1 cm} Since $B^{m}(\chi_\Omega)$ is in
$\operatorname{Lip}(\ep,\Omega)$ for each non-negative number
$m$\,, the Sublemma immediately provides the right control on the
term $II,$ namely,

\begin{equation*}
\begin{split}
II & \le \frac{C}{|\zeta-z_2|^{n+2-l}} \left(|z_1-z_2|^\varepsilon
\, |\zeta-z_2|^{n-l} + |z_1-z_2|^{1+\varepsilon} \,
|\zeta-z_2|^{n-l-1} \right) \\*[7pt] & \le C \,
\frac{|z_1-z_2|^\varepsilon}{|\zeta-z_2|^2}\,,
\end{split}
\end{equation*}
and this completes the proof of Lemma 6.
\end{proof}

\begin{proof}[Proof of the Sublemma]
The most convenient way of proving the Sublemma is to place
ourselves in a real variables context. Given a smooth function $f$
on $\mathbb{R}^d$ let
$$
T_m(f,a)(x) = \sum_{|\alpha| \le m} \frac{\partial^\alpha f
(a)}{\alpha !} \, (x-a)^\alpha
$$
be its Taylor polynomial of degree $m$ around the point $a.$ Then,
clearly,
$$
R_{n-l}(z_1,\zeta)-R_{n-l}(z_2,\zeta) =
T_{n-l}(\frac{\partial^{\ell}h}{\partial \zeta^{\ell}},z_2)(\zeta) -
T_{n-l}(\frac{\partial^{\ell}h}{\partial
\zeta^{\ell}},z_1)(\zeta)\,,
$$
and so the sublemma is an easy consequence of the fact that each
$n$-th order derivative of $h$ is a constant times
$B^m(\chi_\Omega)$\,, for an appropriate exponent $m$\,, and the
following elementary calculus lemma.
\end{proof}

\begin{lemma}\label{calculus}
If $f$ is a $m$ times continuously differentiable function on
$\mathbb{R}^d$, then
\begin{equation*}
\begin{split}
T_m(f,a_1)(x) - T_m(f,a_2)(x) & = \sum_{|\alpha|=m} \partial^\alpha
f (a_1) - \partial^\alpha f (a_2)
\\ & - \sum_{|\alpha|< m} \frac{1}{\alpha!} \left(\partial^\alpha f (a_2)-
T_{m-|\alpha|}\,(\partial^\alpha f,a_1)(a_2)\right) (x-a_2)^\alpha
\,.
\end{split}
\end{equation*}
\end{lemma}
\begin{proof}
Let $P(x)$ stand for the polynomial $T_m(f,a_1)(x) -
T_m(f,a_2)(x),$\, so that
$$
P(x)= \sum_{|\alpha| \le m} \frac{\partial^\alpha P (a_2)}{\alpha !}
\, (x-a_2)^\alpha \,.
$$
A straightforward computation yields
$$
\partial^\alpha P (a_2)= T_{m-|\alpha|}(\partial^\alpha f,a_1)(a_2)-
\partial^\alpha f(a_2)\,,
$$
which completes the proof of Lemma 7.
\end{proof}

\section{\boldmath$\phi$ is bilipschitz}\label{sec6}

In the precedings sections we have proved that $\phi$ is a Lipschitz
function on~$\mathbb{C}$. Moreover
$$
\overline{\partial}\phi=h=(I-\mu B)^{-1} (\mu)\in
\operatorname{Lip}(\varepsilon',\Omega),\quad
0<\varepsilon'<\varepsilon,
$$
and so, by the Main Lemma,
$$
\partial\phi=1+B(h)\in\operatorname{Lip}(\varepsilon',\Omega)\cap \operatorname{Lip}(\varepsilon',(\overline{\Omega})^c), \quad 0<\varepsilon' <\varepsilon.
$$
Since $\phi$ is holomorphic on $(\overline{\Omega})^c$,
$\phi'(z)=\partial \phi (z)$ extends continuosly to $\Omega^c$, but
nothing excludes that this extension might vanish somewhere on
$\partial\Omega$. The functions~$\overline{\partial}\phi$ and
$\partial\phi$ also extend continuously from~$\Omega$ to
$\overline{\Omega}$, but again it could well happen that both vanish
at some point of~$\partial\Omega$. We will show now that this is not
possible. Indeed, we claim that for some positive
number~$\varepsilon_{0}$ we have
\begin{equation}\label{eq26}
|\partial \phi(z)|\ge \varepsilon_{0},\quad z\in\Omega\cap
(\overline{\Omega})^c.
\end{equation}
This implies that the Jacobian of~$\phi$ is bounded from below by
$(1-\|\mu\|^2_{\infty})\,\varepsilon_{0}$ at~$z$ almost all points
of~$\mathbb{C}$. Thus the invers mapping~$\phi^{-1}$ has gradient in
$L^{\infty}(\mathbb{C})$ and hence $\phi$ is bilipschitz.

\begin{proof}[Proof of \eqref{eq26}]
For $a\in\partial\Omega$ denote by $\phi'(a)$ the limit of
$\phi'(z)$ as $z\in (\overline{\Omega})^c$ tends to~$a$. We claim
that \eqref{eq26} follows if we can show that
\begin{equation}\label{eq27}
\phi'(a)\ne 0,\quad a\in\partial\Omega.
\end{equation}
Indeed, this clearly implies
$\inf\limits_{z\in(\overline{\Omega})^c}|\partial\phi(z)|>0$. Now
denote by $\partial\phi(a)$ and $\overline{\partial}\phi(a)$, $a\in
\partial\Omega$, the limits of $\partial\phi(z)$ and
$\overline{\partial}\phi(z)$ as $z\in\Omega$ tends to~$a$. Take a
parametrization~$z(t)$ of $\partial\Omega$ of class~$C^1$, such that
$z'(t)\ne 0$ for all~$t$. Computing $\frac{d}{dt}\phi(z(t))$ in two
different ways,
$$
\phi'(z(t))\, z'(t)=\frac{d}{dt} \phi
(z(t))=\frac{\partial\phi}{\partial z}(z(t))\,
z'(t)+\frac{\partial\phi}{\partial\overline{z}} (z(t))\,
\overline{z'(t)}
$$
and so
$$
\phi'(z(t))=\left(1+\mu (z(t))\frac{\overline{z'(t)}}{z'(t)}\right)
\frac{\partial\phi}{\partial z}(z(t)).
$$
Thus, by \eqref{eq27}, $\frac{\partial\phi}{\partial z}(a)\ne 0$,
$a\in\partial\Omega$, which yields
$\inf\limits_{z\in\Omega}\left|\frac{\partial\phi}{\partial
z}(z)\right|>0$.

We turn now to the proof of~\eqref{eq27}. Assume that
$0=a\in\partial\Omega$. Performing a rotation before applying $\phi$
we may assume that $\lambda=\mu(0)$ is a non-negative real number
($\mu(0)$ is the limit of~$\mu(z)$ as $z\in\Omega$ tends to~$0$).
Performing a rotation after applying $\phi$ we may also assume that
the tangent plane to~$\partial\Omega$ at the origin is the real
axis. Denote by~$H^+$ and $H^-$ the upper and lower half planes,
respectively. Consider the continuous piecewise linear mapping
$$
z=L(w)=(w-\lambda
\overline{w})\,\chi_{H^-}(w)+(1-\lambda)w\,\chi_{H^+}(w).
$$
Then, by \cite[(5.6), p.~83]{LV}, the Beltrami coeficient~$\nu(w)$
of the mapping~$\phi\circ L$ is
\begin{equation}\label{eq28}
\nu(w)=\frac{\overline{\partial}(\phi\circ L)(w)}{\partial
(\phi\circ L)(w)}=\frac{-\lambda \,\chi_{H^-}(w)+\mu
(L(w))}{1-\lambda\, \chi_{H^-}(w)\,\mu (L(w))}.
\end{equation}
Since $\nu$ vanishes on $H^+\cap L^{-1}((\overline{\Omega})^c)$,
\begin{equation*}
\begin{split}
\int_{|w|<r_{0}}\frac{|\nu(w)|}{|w|^2}\,dA(w)&=\int_{H^+\cap
L^{-1}(\Omega)\cap B(0,r_{0})} \dotsi\\*[7pt] &\quad +\int_{H^-\cap
L^{-1}(\Omega^c)\cap B(0,r_{0})} \dotsi\\*[7pt] &\quad
-\int_{H^-\cap L^{-1}(\Omega)\cap B(0,r_{0})} \dotsi\\*[7pt]
&=I+II+III,
\end{split}
\end{equation*}
where $r_{0}$ is the small number introduced in the proof of the
Main Lemma (see Figure~1). By \eqref{eq28},
$$
|I|\le \int_{H^+\cap L^{-1} (\Omega)\cap B(0,r_{0})}\frac{|\mu
(L(w))|}{|w|^2}\,dA(w).
$$
Since $L(w)=(1-\lambda)w$ on $H^+$, making the change of variables
$z=L(w)$ gives
$$
|I|\le \int_{H^+\cap\Omega \cap B(0,r_{0})}\frac{dA(z)}{|z|^2}\le
\int_{0}^{r_{0}}r^{-1+\varepsilon}\,dr <\infty,
$$
where in the next to the last inequality we used \eqref{eq12}.

For $II$ we begin by remarking that
$$
|II|=\lambda\int_{H^-\cap L^{-1}(\Omega^c)\cap
B(0,r_{0})}\frac{dA(w)}{|w|^2},
$$
and making the change of variables~$z=w-\lambda\,\overline{w}$, we
get
$$
|II|\le\lambda\int_{H\cap \Omega^c\cap
B(0,r_{0})}\frac{1}{|z|^2}\frac{1+\lambda}{1-\lambda}\,dA(z),
$$
which can be shown to be finite as before (in particular, using
again \eqref{eq12}). To take care of~$III$ we  make the same change
of variables and we obtain
\begin{equation*}
\begin{split}
|III|&\le \int_{H^-\cap L^{-1}(\Omega)\cap B(0,r_{0})}
\left|\frac{-\lambda+\mu (L(w))}{1-\lambda\, \mu (L(w))}\right|
\frac{dA(w)}{|w|^2}\\*[7pt] &\le \frac{1}{1-\lambda} \int_{H^-\cap
L^{-1}(\Omega)\cap B(0,r_{0})} \frac{|\mu
(L(w))-\mu(0)|}{|w|^2}\,dA(w)\\*[7pt] &\le
\frac{1+\lambda}{(1-\lambda)^2} \int_{H^-\cap \Omega\cap B(0,r_{0})}
\frac{|\mu (z) -\mu(0)|}{|z|^2}\,dA(z)\\*[7pt] &\le
C\int_{B(0,r_{0})}\frac{dA(z)}{|z|^{2-\varepsilon}}<\infty.
\end{split}
\end{equation*}
Therefore
$$
\int_{|w|<r_{0}}\frac{|\nu(w)|}{|w|^2}\,dA(w)<\infty,
$$
and so, by \cite[p.~232]{LV}, $H=\phi\circ L$ is conformal at the
origin, in the sense that the limit
$$
H'(0)=\lim_{z\to 0}\frac{H(z)-H(0)}{z}
$$
exists and $H'(0)\ne 0$.

The part of the imaginary positive axis close to the origin is
included in~$(\overline{\Omega})^c$ (see Figure~1), and thus
$L^{-1}(iy)=\frac{iy}{1-\lambda}$ if $y>0$ is small. Hence
\begin{equation*}
\begin{split}
\phi'(0)&=\lim_{0<y\to 0}\phi'(iy)=\lim_{0<y\to 0}
\frac{\phi(iy)-\phi(0)}{iy}\\*[7pt] &=\lim_{0<y\to
0}\frac{H(L^{-1}(iy))-H(0)}{iy}\\*[7pt] &=
\frac{H'(0)}{(1-\lambda)}.
\end{split}
\end{equation*}
This completes the proof of \eqref{eq27}.
\end{proof}

\section{Reduction to the one domain case}
Suppose, as in the statement of the Theorem, that
$\Omega_1,$\,...\,,$\Omega_N$ are bounded disjoint domains with
boundary of class ~$C^{1+\ep}$, for some $\ep$ with $0<\ep <1$,
and that $\mu = \sum_{j=1}^N \mu_j\,\chi_{\Omega_j},$ where
$\mu_j$ is of class $\operatorname{Lip}(\ep,\Omega_j),$\, and
$\|\mu\|_\infty < 1 \,.$ Let $\Phi^\mu$ be the quasiconformal
mapping associated with $\mu\,.$ Assume that $N > 1$ and set
$\nu_1 = \mu_1 \,\chi_{\Omega_1}$ and $\nu_2 = \sum_{j=2}^N
\mu_j\,\chi_{\Omega_j}\,.$  By \cite[(10) p. 9]{Ah},
$$\Phi^\mu = \Phi^\lambda \circ \Phi^{\nu_1}\,,$$
 where
$$
\lambda(w)= \nu_2(z) \, \frac{\partial \Phi
^{\nu_1}(z)}{\overline{\partial \Phi ^{\nu_1}(z)}}\,,
$$
and, for each $w \in \C$ the point $z$ is defined by $ w=
\Phi^{\nu_1}(z)\,.$ In particular, $\lambda$ is supported on $
\cup_{j=2}^N \Phi^{\nu_1}(\Omega_j)\,.$ Recall from the previous
section that $\Phi^{\nu_1}$ is bilipschitz, of class
$C^{1+\varepsilon'}\,,\, 0 <\varepsilon' <\varepsilon, $ on
$\Omega_1$ and $(\overline{\Omega_1})^c \,,$\, and conformal on
$\Omega_1^c \,.$  In particular, $\Phi^{\nu_1}$ is holomorphic on
$(\overline{\Omega_1})^c$ and
$$\frac{d \Phi^{\nu_1}}{d z}(z) \neq
0\,, \,\, z \in \Omega_1^c \,.$$ Thus the bounded domains
$\Phi^{\nu_1}(\Omega_j)\,,\, 1 \le j \le N \,,$  have boundaries
of class $C^{1+\varepsilon'}\,,\, 0 < \varepsilon' < \varepsilon
\,.$ On the other hand, the function $\lambda$ satisfies a
Lipschitz condition of order $\varepsilon'\,,$\, for $ 0 <
\varepsilon' < \varepsilon \,,$ in each domain
$\Phi^{\nu_1}(\Omega_j)\,,\, 2 \le j \le N \,.$ Proceeding by
induction we conclude now that $\Phi^\mu$ is bilipschitz and of
class $C^{1+\varepsilon'}\,,\, 0 < \varepsilon' < \varepsilon \,,$
in each domain $\Omega_j\,,\,  1 \le j \le N \,.$

One can prove now Corollary 2. Let $D$ be a bounded planar domain
and let $f$ be a function in $W^{1,2}_{\text{loc}}(D)$ satisfying
the Beltrami equation
\begin{equation*}
\frac{\partial f}{\partial \overline{z}}(z)= \mu(z)\, \frac{\partial
f}{\partial z}(z)\,,\quad z \in D\,,
\end{equation*}
where $\mu$ is as in the statement of the Theorem. By Stoilow's
factorization Theorem, $ f= h \circ \Phi\,,$ where $\Phi$ is the
quasiconformal mapping associated with $\mu$ and $h$ is a
holomorphic function on $\Phi(D)\,.$  Let $D_\delta$ stand for the
set of points in $D$ whose distance to the boundary of $D$ is larger
than $\delta\,.$ Thus $h$  has bounded derivatives of all orders on
$\Phi(D_\delta)$ and consequently $f$ is as smooth as $\Phi$ on
$D_\delta\,.$  Hence $f$ is Lipschitz on $D_\delta$ and of class
$C^{1+\varepsilon'}\,,\, 0 < \varepsilon' < \varepsilon \,,$ in each
open set $D_\delta \cap \Omega_j\,,\,  1 \le j \le N \,.$ If $D$
contains the closure of each $\Omega_j\,,$ then $f$ is of class
$C^{1+\varepsilon'}\,,\, 0 < \varepsilon' < \varepsilon \,,$ on each
$\Omega_j \,.$

Recalling the relation between the Beltrami equation and second
order elliptic equations in divergence form, as explained in the
introduction, Corollary 2 follows immediately from the above
argument.

\section{Cuspidal domains}
As we remarked in the introduction, the conclusion of the Theorem
fails for domains with corners; for instance, for a square. However,
the class of domains with boundary of class $C^{1+\ep}$ is not
optimal for the Theorem. There is a heuristic argument that points
out at a more general class of domains, which, at least in a first
approximation,  may be viewed as optimal.

First of all we recall that a central point in the proof of the
Theorem was the fact, which is part of the Main Lemma, that each
power of the Beurling transform sends the characteristic function
of the domain into a bounded function. Let us concentrate on the
Beurling transform $B$ and find a simple condition on a bounded
domain $\Omega$ with rectifiable boundary so that $B(\chi_\Omega)$
is bounded. Our first remark is that $B(\chi_\Omega)$ can be
written as the Cauchy transform of a boundary measure. For this we
use, on one hand, the basic property of $B\,,$
$$
\partial \varphi = B(\overline{\partial} \varphi)\,,
$$
which holds for all compactly supported smooth functions $\varphi$
and extends to a variety of situations by regularization. On the
other hand, we use the elementary identity
$$
\partial \chi_\Omega = \frac{1}{2\,\imath}\,\, d
\overline{z}_{\partial\Omega} \,,
$$
which holds at least for bounded domains with rectifiable boundary.
Combining the above two identities we get
$$
\overline{\partial}\,B(\chi_\Omega)=
B(\overline{\partial}\chi_\Omega) = \partial \chi_\Omega =
\frac{1}{2\,\imath}\, d \overline{z}_{\partial\Omega} \,,
$$
which yields
$$
B(\chi_\Omega) = \frac{1}{2\,\imath}\, C(d
\overline{z}_{\partial\Omega})\,.
$$
Now, $d \overline{z}_{\partial\Omega} = \overline{\tau}^2 (z)\,d
z_{\partial\Omega}\,,$\, $\tau(z)$ being the unit tangent vector
to $\partial\Omega$ at $z\,.$  Assume now that the arc-length
measure on the boundary of $\Omega$ satisfies the Ahlfors
condition $\text{length}(\partial\Omega \cap D(z,r)) \le C\, r\,,$
for each $z \in \partial\Omega$ and $r>0\,,$\, where $D(z,r)$
stands for the open disc with center $z$ and radius $r\,.$ Then a
simple estimate shows that the Cauchy integral of a function $f$
on $\partial\Omega\,,$ that is,
$$\frac{1}{2\pi \imath}\,\int_{\partial\Omega} \frac{f(w)}{z-w}\,dw \,,
\quad z \in \C \setminus \partial\Omega\,,$$ is bounded provided
$f$ satisfies a Lipschitz condition of some positive order on
$\partial\Omega\,.$  Thehefore, to get boundedness of
$B(\chi_\Omega)$ one has to require that the square of the tangent
unit vector satisfies a Lipschitz condition of some positive order
on $\partial\Omega\,.$  This is weaker than requiring the
Lipschitz condition on the tangent unit vector itself, because it
allows jumps of $180$ degrees on the argument of the tangent unit
vector. In other words, cusps are allowed.

We want now to define formally \textit{cuspidal domains of class}
$C^{1+\ep}\,.$  Given a planar domain $\Omega$ we say that
$\partial\Omega$ is $C^{1+\ep}$-smooth at a boundary point $z_0$
if there is a positive $r_0$ such that $\Omega\cap D(z_0,r_0)$ is,
after possibly a rotation, the part of $D(z_0,r_0)$ lying below
the graph of a function of class $C^{1+\ep}\,.$

 We say that $\Omega$ has an
\textit{interior cusp of class} $C^{1+\ep}$ at $z_0= x_0+\imath
y_0\in
\partial\Omega$ provided there is a positive $r_0$  and functions
$y=a(x)\,,$ $y= b(x)\,,$\, of class $C^{1+\ep}$ on the interval
$(x_0-r_0,x_0+r_0)\,,$ such that $a(x_0)=
a'(x_0)=b(x_0)=b'(x_0)=0$, and, after possibly a rotation, a point
$z=x+\imath y$ is in $\Omega\cap D(z_0,r_0)$ if and only if $z \in
D(z_0,r_0)$ and $b(x) < y  < a(x)\,.$

We say that $\Omega$ has an \textit{exterior cusp of class}
$C^{1+\ep}$ at $z_0 \in \partial\Omega$  provided
$\overline{\Omega}^c$ has an interior cusp of class $C^{1+\ep}$ at
$z_0\,.$

A planar domain $\Omega$ is a \textit{cuspidal domain of class}
$C^{1+\ep}$ if $\partial\Omega$ is $C^{1+\ep}$-smooth at all
boundary points, except possibly at finitely many boundary points
where $\Omega$ has a cusp of class $C^{1+\ep}$ (either interior or
exterior).

The simplest examples of non-smooth cuspidal domains are the drop
like domain $D$ and the peach like domain $P$ shown in figure 5
below. The reader may easily imagine more complicated cuspidal
domains with lots of cusps of both types (see figure 7).

\begin{figure}[ht]
\begin{center}
\includegraphics{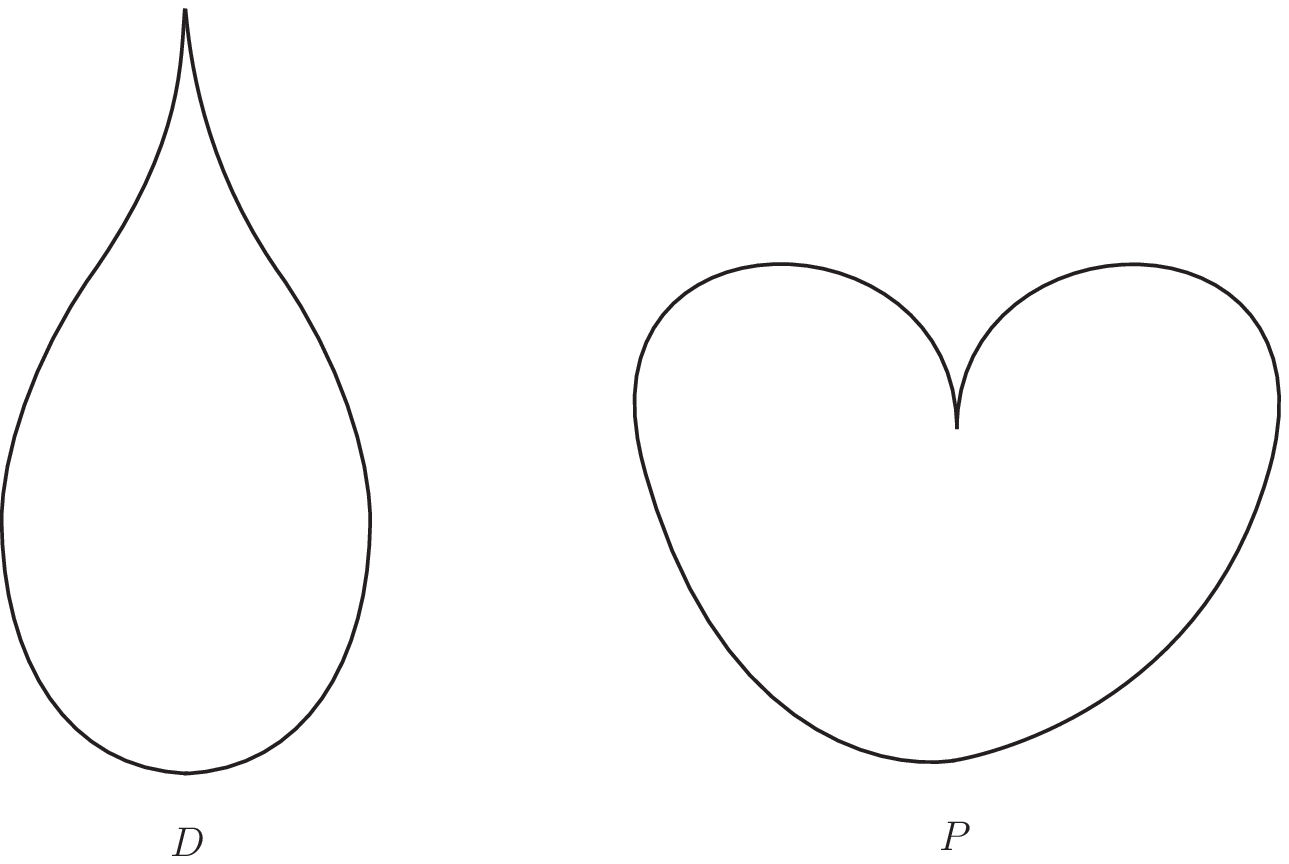}\\*[7pt]
Figure 5
\end{center}
\end{figure}

With appropriate formulations the Theorem and Corollaries 1 and 2
hold true for cuspidal domains of class $C^{1+\ep}\,.$ The right
statements involve the notion of geodesic distance in the domain
$\Omega$, which we discuss now. Given two points $z$ and $w$ in
$\Omega$ their geodesic distance is defined by
$$
d(z,w)= d_\Omega(z,w)= \inf_{\gamma} l(\gamma)\,,
$$
where the infimum is taken over all rectifiable curves $\gamma$ in
$\Omega$ joining $z$ and $w\,.$ Here $l(\gamma)$ stands for the
length of $\gamma\,.$  Notice that if $\Omega$ has only interior
cusps, then the geodesic and the Euclidean distances are
comparable and that this is not the case in the proximity of an
exterior cusp. The Lipschitz norm
$\|\cdot\|_{\varepsilon,\,\Omega,\,d_\Omega}$ of order
$\varepsilon$ and the corresponding Lipschitz spaces
$\operatorname{Lip}(\varepsilon,\Omega,d_\Omega)$ with respect to
the distance $d$ are defined is the usual way, with the Euclidean
distance replaced by $d$ in \eqref{eq2} and \eqref{eq2bis}.

The Theorem for cuspidal domains reads as follows.
\begin{teor'}
Let $\{\Omega_j\}$\,, $1 \le j \le N$\,, be a finite family of
disjoint bounded cuspidal domains of class~$C^{1+\ep}$, $0<\ep
<1$, and let $\mu = \sum_{j=1}^N \mu_j\,\chi_{\Omega_j},$ where
$\mu_j$ is of class
$\operatorname{Lip}(\ep,\Omega_j,d_{\Omega_j}).$ \, Assume in
addition that $\|\mu\|_\infty < 1$. Then the associated
quasiconformal mapping $\Phi$ is bilipschitz.
\end{teor'}
Corollary 1 remains true without any change.
\begin{cor1} \label{C1'}
If $\Omega$ is a bounded cuspidal domain of class $C^{1+\ep}$,
$0<\ep <1$, and $\mu= \lambda \,\chi_\Omega$, where $\lambda$ is a
complex number such that $|\lambda|< 1$, then the associated
quasiconformal mapping~$\Phi$ is bilipschitz.
\end{cor1}

\begin{cor2}\label{C2'}
Let $\Omega_j$\,, $1 \le j \le N$\,, be a finite family of
disjoint bounded cuspidal domains of class~$C^{1+\ep}$, $0<\ep
<1$, and assume that all domains $\Omega_j$ are contained in a
bounded domain $D$ with boundary of class~$C^{1+\ep}$. Let
$A=A(z),\,\, z \in D,$\, a $2\times 2$ symmetric elliptic matrix
with determinant $1$ and entries supported in $ \cup_{j=1}^N
\Omega_j$ and belonging to
$\operatorname{Lip}(\ep,\Omega_j,d_{\Omega_j}),$\,\, $1 \le j \le
N .$ Let $u$ be a solution of equation \eqref{eq1bis} in $D.$ Let
$D_\delta$ stand for the set of points in $D$ at distance greater
than $\delta$ from the boundary of $D.$ Then $\nabla u \in
\operatorname{Lip}(\varepsilon',\Omega_j \cap D_\delta
,d_{\Omega_j}) ,$\,\, for $0 <\varepsilon' < \varepsilon $ and
\,$1 \le j \le N .$ In particular, $\nabla u \in
L^\infty(D_\delta)$ and $u$ is a locally Lipschitz function in
$D.$
\end{cor2}

 We proceed now to sketch the proof of Theorem'. The modifications needed
 are minor and fortunately one can  reduce without much pain the cuspidal
case to the smooth case.

We start by discussing the proof in the one domain case ($N=1$).
The difficulty is that the Main Lemma does not hold for cuspidal
domains with exterior cusps and the usual Euclidean Lipschitz
spaces. We present an example to make the difficulty clear.

\begin{example}
Let $\Omega$ be the open disc centered at the origin of radius $2$
minus the union of the closed discs of radius $1$ centered at $1$
and $-1\,.$ Thus $\Omega$ has an exterior cusp at the origin. We
claim that $B(\chi_\Omega)$ does not satisfy a Lipschitz condition
on $\Omega$ of any order $\varepsilon$ such that $1/2 <\varepsilon
\leq 1 \,.$ The point is that we can explicitly calculate
$B(\chi_\Omega)\,.$ If $D(a,r)$ stands for the open disc centered
at $a$ of radius $r\,,$ then a simple argument (see, for example, \cite[p.~965]{MV})
shows that
$$
B(\chi_{D(a,r)})(z) =
-\frac{r^2}{(z-a)^2}\,\chi_{D^c(a,r)}(z)\,,\quad z \in \C\,.
$$
Thus
$$
B(\chi_\Omega)(z)= \frac{1}{(z-1)^2}+ \frac{1}{(z+1)^2}\,, \quad z
\in \C \setminus (\overline{D(-1,1)}\cup \overline{D(1,1)})\,.
$$
Set $z_1 = -x +\imath y$ and $z_2 = x +\imath y$\,, where $x$ and
$y$ are positive real numbers such that $z_1$ and $z_2$ are in
$\Omega$\,. Then $|z_1-z_2|= 2\,x \,.$ On the other hand, a simple
computation yields $|B(\chi_\Omega)(z_1)-B(\chi_\Omega)(z_2)|
\simeq y$ as $x$ and $y$ tend to $0\,.$ Choosing $ y \simeq
\sqrt{x}$ we conclude that $B(\chi_\Omega)$ does not satisfy a
Lipschitz condition on $\Omega$ of any order $\varepsilon$ with
$1/2 <\varepsilon \leq 1 \,. $
\end{example}

The Main Lemma for cuspidal domains reads as follows.
\begin{ML}
Let $\Omega$ be a bounded cuspidal domain  of class $C^{1+\ep}$,
$0<\ep<1$, and let $T$ be an even smooth homogeneous
Calder\'{o}n-Zygmund operator. Then $T$ maps
$\operatorname{Lip}(\varepsilon,\Omega,d_\Omega)$ into itself, and
 $T$ also maps $\operatorname{Lip}(\varepsilon,\Omega,d_\Omega)$ into
  $\operatorname{Lip}(\ep,\overline{\Omega}^c, d_{\overline{\Omega}^c})$. In fact, one
has the inequalities
$$
\| Tf \|_{\ep,\,\Omega,\,d_\Omega} \leq C\,
\|T\|_{CZ}\,\|f\|_{\ep,\,\Omega,\,d_\Omega}
$$
and
$$
\| Tf \|_{\ep,\,\overline{\Omega}^c,\,d_{\overline{\Omega}^c}}
\leq C\, \|T\|_{CZ}\,\|f\|_{\ep,\,\Omega,\,d_\Omega}
$$
 where $C$ is a constant depending only on $\ep$
and $\Omega$\,.
\end{ML}
\begin{proof}
We first prove that $T^*(f) \in L^\infty(\C)$ for each $f \in
\operatorname{Lip}(\varepsilon,\Omega,d_\Omega)\,.$ Take a point
$z \in \C.$  Clearly we may assume that $z \in  D(z_0,\frac{1}{3}
\,r_0)\,,$ where $z_0$ is a cuspidal point and $r_0$ is as in the
definition of cusp. Otherwise $z$ is far from all cusps and thus
we may apply the arguments of the smooth case. Since there are
only finitely many cusps we may also assume that there is a
positive number $r_0$ which works in the definition of cusp
simultaneously for all cusps. Clearly
$$
T^*(f)(z) \le T^*(f\,\chi_{\Omega \cap D(z_0,r_0)})(z) +
T^*(f\,\chi_{\Omega \cap D^c(z_0,r_0)})(z)\,,
$$
and the second term is bounded by
$C\,\log{\frac{d}{r_0}}\,\|f\|_\infty\,,$ where $d$ stands for the
diameter of $\Omega.$  We are therefore left with the first term.

Assume for the moment that $\Omega$ has an interior cusp at
$z_0\,.$ We connect the cercle of center $z_0$ and radius
$\frac{2}{3}\,r_0$ with the concentric cercle of radius $r_0$ to
produce domains $D_1$ and $D_2$ with boundary of class
$C^{1+\varepsilon}\,,$ as shown in figure 6.

\begin{figure}[ht]
\begin{center}
\includegraphics{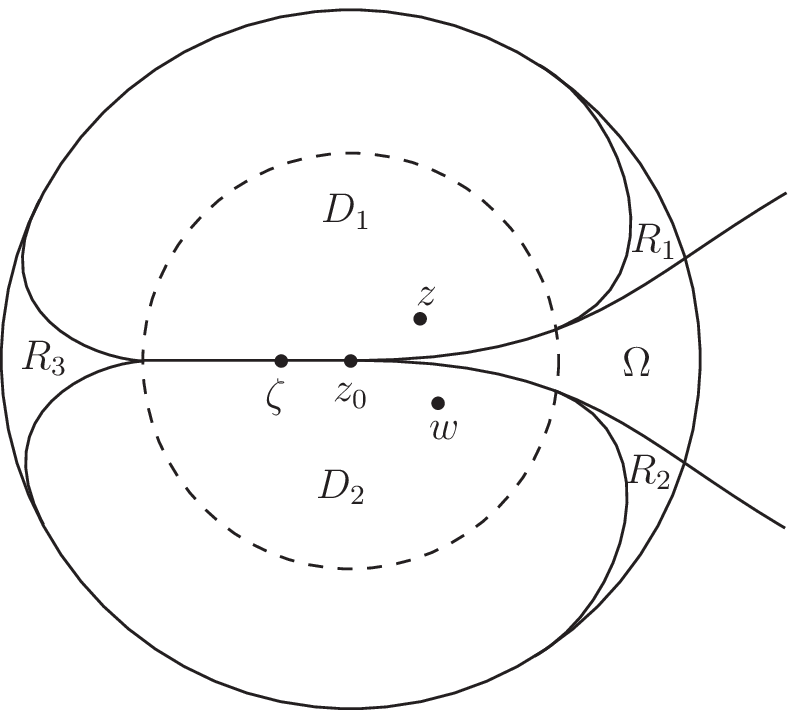}\\*[7pt]
Figure 6
\end{center}
\end{figure}

\noindent
Notice that three
other ``residual" domains $R_j\,,$\, $1 \leq j \leq 3\,,$ have
been formed. The domains $R_j$ are not smooth but they are far
from $D(z_0,\frac{1}{3}\,r_0)\,.$ Since we are assuming that
$\Omega$ has an interior cusp at $z_0\,,$ the restriction of $f$
to $\Omega \cap D(z_0,r_0)$ satisfies an Euclidean Lipschitz
condition of order $\varepsilon\,.$  By the well known extension
theorem for Lipschitz functions (\cite[Chapter VI ]{St}) we may
extend the restriction of $f$ to  $\Omega \cap D(z_0,r_0)$ to a
function $g \in \operatorname{Lip}(\varepsilon,\C)$ such that
\begin{equation}\label{eq31}
\|g\|_{\varepsilon,\,\C} \le C\, \|f\|_{\varepsilon,\,\Omega\cap
D(z_0,r_0)}\,.
\end{equation}
Hence
$$
f\,\chi_{\Omega \cap D(z_0,r_0)} = g\,\chi_{D(z_0,r_0)} -
g\,\chi_{D_1}-g\,\chi_{D_2} - \sum_{j=1}^3 g\,\chi_{R_j}
$$
and so
$$
T^*(f\,\chi_{\Omega \cap D(z_0,r_0)}) \le
T^*(g\,\chi_{D(z_0,r_0)})+ T^*(g\,\chi_{D_1})+
T^*(g\,\chi_{D_2})+\sum_{j=1}^3 T^*(g\,\chi_{R_j})\,.
$$
The last three terms are controlled by $C\,\|g\|_\infty\,,$ where
$C$ is a constant depending on $r_0\,,$ because the domains of the
functions to which $T^*$ is applied are far from $D(z_0,
\frac{1}{3}\,r_0)\,.$ By \eqref{eq31} this gives the correct bound
$C\,\|f \|_{\varepsilon,\Omega} \,.$

To estimate the first three terms we remark that $D_1\,,$ $D_2$
and $D(z_0,r_0)$ are domains with boundary of class
$C^{1+\varepsilon}\,.$ Thus the proof of the Main Lemma for smooth
domains of section 3 yields, by \eqref{eq31}\,,
\begin{equation}\label{eq32}
\|T^*(f)\|_\infty \le C\, \|f \|_{\varepsilon,\Omega}\,.
\end{equation}

Assume now that $\Omega$ has an exterior cusp at $z_0\,.$ Then
$$
f\,\chi_{\Omega \cap D(z_0,r_0)} = f\,\chi_{D_1}+ f \,\chi_{D_2} +
\sum_{j=1}^3 f \,\chi_{R_j}
$$
and so
$$
T^*(f\,\chi_{\Omega \cap D(z_0,r_0)}) \le  T^*(f\,\chi_{D_1})+
T^*(f\,\chi_{D_2})+\sum_{j=1}^3 T^*(f\,\chi_{R_j})\,.
$$
The proof of the estimate \eqref{eq32} proceeds as before.

Let us turn to prove the Lipschitz condition on $T(f)$. Since
$T(f)$ is bounded, to estimate $|T(f)(z)-T(f)(w)|$ we may restrict
our attention to the case in which $z$ and $w$ are very close to
each other. We may also assume that $z$ and $w$ are close to a
cusp, say, $z, w \in D(z_0,\frac{1}{3}\,r_0 )\,.$  From this point
on the proof is very similar to what we did before, with $T^*$
replaced by $T\,.$ The first step is to restrict our attention to
$f\,\chi_{\Omega\cap D(z_0,r_0)}\,,$ which may be achieved by the
identity
$$
T(f) = T(f\,\chi_{\Omega \cap D(z_0,r_0)}) + T(f\,\chi_{\Omega
\cap D^c(z_0,r_0)})\,.
$$

Assume first that $\Omega$ has an interior cusp at $z_0$ and
consider again the extension $g$ of $f\,\chi_{\Omega\cap
D(z_0,r_0)}$ satisfying \eqref{eq31}. We clearly have
$$
T(f\,\chi_{\Omega \cap D(z_0,r_0)})=
T(g\,\chi_{D(z_0,r_0)})-T(g\,\chi_{D_1})-T(g\,\chi_{D_2})-\sum_{j=1}^3
T(g\,\chi_{R_j})\,.
$$
By the smooth version of the Main Lemma $T(g\,\chi_{D_1})$ and
$T(g\,\chi_{D_2})$ satisfy an Euclidean Lipschitz condition of order
$\varepsilon$ on the complement of ${D_1}$ and ${D_2}$ respectively.
Clearly, again by the smooth Main Lemma, $T(g\,\chi_{D(z_0,r_0)})$
satisfies an Euclidean Lipschitz condition of order $\varepsilon$ on
$D(z_0,r_0)\,.$ Finally $T(g\,\chi_{R_j})\,, 1 \le j\le 3\,,$
satisfy an Euclidean Lipschitz condition of order $\varepsilon$ on
$D(z_0,\frac{1}{3}r_0)\,,$ because the domains $R_j$ are far from
$D(z_0,\frac{1}{3}r_0)\,.$ If both $z$ and $w$ belong to $\Omega$ we
then get an estimate of the form
$$
|T(f)(z)-T(f)(w)| \le
C\,\|f\|_{\varepsilon,\Omega}\,|z-w|^\varepsilon\,,
$$
which completes the proof, because in the case at hand the
Euclidean distance between $z$ and $w$ is comparable to their
geodesic distance. If $z$ and $w$ are in $\overline{\Omega}^c\,,$
then we may assume that $z \in D_1$ and $w \in D_2\,.$ Otherwise
$z$ and $w$ belong both to either $D_1$ or $D_2$ and thus we
obtain the above Euclidean Lipschitz estimate of order
$\varepsilon$. Choose then a point $\xi \in
\partial{D_1} \cap
\partial{D_2}$ such that $\max \{|z-\xi|,|w-\xi| \} \simeq
d_{\overline{\Omega}^c}(z,w)\,.$ (see figure 6). Since
$T(g\,\chi_{D_1})$ and $T(g\,\chi_{D_2})$ satisfy an Euclidean
Lipschitz condition of order $\varepsilon$ on $D_1$ and $D_2$
respectively, and since $T(f)$ is continuous on the complement of
$\Omega\,,$ because $f$ is supported on $\Omega\,,$ we get
$$
|T(f)(z)-T(f)(w)| \le C\,\max
\{|z-\xi|^\varepsilon,|w-\xi|^\varepsilon\}\simeq
d_{\overline{\Omega}^c}(z,w)^\varepsilon\,.
$$

If $\Omega$ has an exterior cusp at $z_0\,,$ we argue similarly,
using the identity
$$
T(f\,\chi_{\Omega \cap D(z_0,r_0)})=
T(f\,\chi_{D_1})+T(f\,\chi_{D_2})+\sum_{j=1}^3 T(f\,\chi_{R_j})\,.
$$
 The details are left to the reader.
\end{proof}

We continue now  the proof of Theorem' in the one domain case. If
the domain has only interior cusps the argument we described to
prove the Theorem goes through without any change. The reason is
that, since the geodesic distance in $\Omega$ is comparable to the
Euclidean distance sections 4 and 5 hold true, because of the
cuspidal version of the Main Lemma. For section 6 one has to
remark that $\nabla \Phi$ is now only in
$\operatorname{Lip}(\ep,\overline{\Omega}^c,
d_{\overline{\Omega}^c})\,,$ but that this still implies that
$\nabla \Phi$ extends continuously from $\overline{\Omega}^c$ to
$\partial\Omega\,.$ The conformality of $\Phi$ on $\Omega^c$ is
proved as in section 6, after remarking that cusps do not create
any problem because they are perfectly suited for an appeal to
\cite[p.232]{LV}.

Assume now that our domain has exterior cusps. In figure 7 it is
shown how to subdivide $\Omega$ in finitely many subdomains
$\Omega_j\,,\, 1 \le j \le M$, which are cuspidal domains of class
$C^{1+\varepsilon}$ \textit{with only interior cusps}. Since $\mu
\in \operatorname{Lip}(\varepsilon,\Omega,d_\Omega)\,,$ $\mu_j =
\mu \,\chi_{\Omega_j} $ is in
$\operatorname{Lip}(\varepsilon,\Omega_j)\,.$ We can then apply
the factorization method of section 7 to get that $\nabla\Phi \in
\operatorname{Lip}(\varepsilon',\Omega,d_\Omega)\,$ and
$\nabla\Phi \in
\operatorname{Lip}(\varepsilon',\overline{\Omega}^c,d_{\overline{\Omega}^c})\,,$
\,$0 < \varepsilon' <\varepsilon \,.$

\begin{figure}[ht]
\begin{center}
\includegraphics{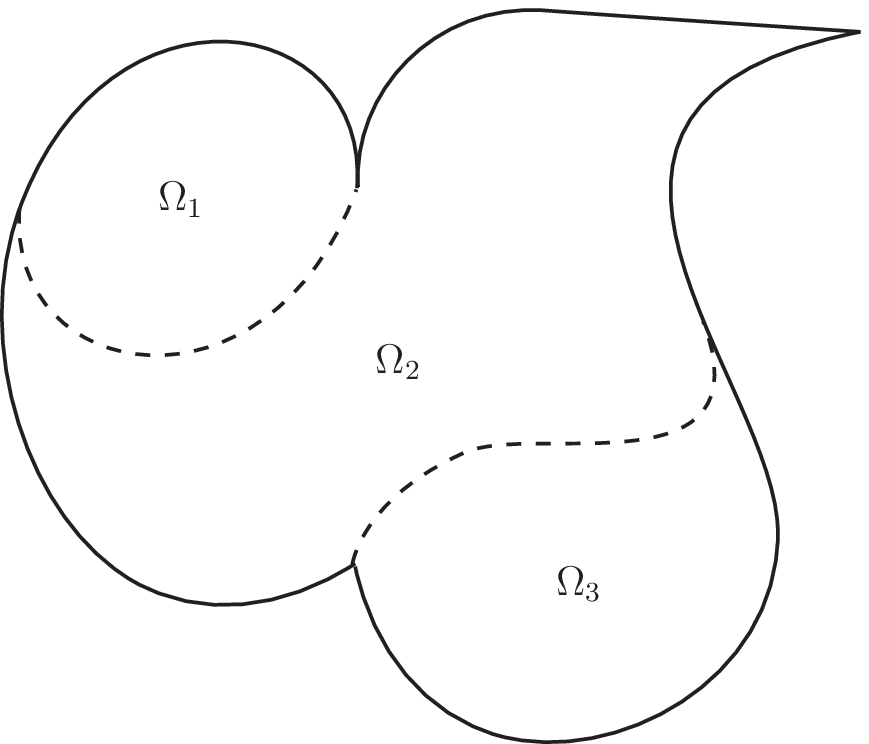}\\*[7pt]
Figure 7
\end{center}
\end{figure}

The reduction to the one domain case needs only one comment. Using
the notation of section 7, we use the fact that $\Phi^{\nu_1}$ is
conformal on $\Omega_1^c$ to ascertain that the image of each
$\Omega_j$ under $\Phi^{\nu_1}$ is again a cuspidal domain of
class
 $C^{1+\varepsilon'}\,,\, 0 < \varepsilon' < \varepsilon \,.$  We repeat for emphasis that
 the conformality of  $\Phi^{\nu_1}$ is proved at a cusp appealing, as
in section 6, to \cite[p.~232]{LV}.

\section{Final comments}
Very likely the restriction on the determinant of the matrix $A$
in Corollary 2 and Corollary 2' is superfluous. This would follow
if Lipschitz regularity results should hold for the general
elliptic system
$$
\overline{\partial}\Phi(z) = \mu(z)\,\partial\Phi(z)+
\nu(z)\overline{\partial\Phi(z)}\,,
$$
where $|\mu(z)|+|\nu(z)| \le k<1\,,\,$ a. e. on $\C\,.$

It seems also rather clear that the right conclusion in Corollary
2 (and analogously in Corollary 2') should be that the solution $u
$ is of class $\operatorname{Lip}(\varepsilon,\Omega_j)\,,$\, $1
\le j \le N \,,$ in all dimensions (and without any restriction on
the determinant of $A$). Evidence for this conjecture is provided
by the fact that it is true in the plane whenever the Lipschitz
norm of the coefficients of the matrix $A$ are small enough. See
the argument for the Beltrami equation at the end of section 2. We
acknowledge some useful correspondence with L. Escauriaza and
D.Faraco on that issue.

Apparently it is not known what is the best exponent $p$ such that
$\nabla\Phi^\mu \in L^p_{\text{loc}}(\C)$ for $\mu = \lambda\,
\chi_Q\,,$\, where $\lambda $ is a complex number such that $|\lambda|<1$ and $Q$ is a square. This looks
extremely surprising to the authors, who would very much
appreciate knowing the precise regularity properties of $\Phi^\mu$
in the scale of the local Sobolev spaces
$W^{1,p}_{\text{loc}}(\C)$.

\begin{gracies}
This work was partially supported by the grants \newline
2005SGR00774 (Generalitat de Catalunya) and  MTM2007-60062
(Ministerio de Educaci\'{o}n y Ciencia). The authors are indebted to
Stephen Semmes for suggesting the proof of \eqref{eq27} and to
Eero Saksman for an encouraging conversation.
\end{gracies}

\begin{tabular}{l}
Joan Mateu\\
Departament de Matem\`{a}tiques\\
Universitat Aut\`{o}noma de Barcelona\\
08193 Bellaterra, Barcelona, Catalonia\\
{\it E-mail:} {\tt mateu@mat.uab.cat}\\ \\
Joan Orobitg\\
Departament de Matem\`{a}tiques\\
Universitat Aut\`{o}noma de Barcelona\\
08193 Bellaterra, Barcelona, Catalonia\\
{\it E-mail:} {\tt orobitg@mat.uab.cat}\\ \\
Joan Verdera\\
Departament de Matem\`{a}tiques\\
Universitat Aut\`{o}noma de Barcelona\\
08193 Bellaterra, Barcelona, Catalonia\\
{\it E-mail:} {\tt jvm@mat.uab.cat}
\end{tabular}
\end{document}